%% file: bol.tex
\documentclass[francais]{smfart}

\usepackage[dvips]{epsfig}

\usepackage[francais]{babel}

\usepackage{smfthm}

\newcommand{\R}{\mathbb R}
\newcommand{\C}{\mathbb C}
\renewcommand{\P}{\mathbb P}
\newcommand{\N}{\mathbb N}
\newcommand{\Z}{\mathbb Z}

\newcommand{\G}{\mathbb G}

\newtheorem{theorem}{Th\'eor\`eme}[section]
\newtheorem{lemma}[theorem]{Lemme}
\newtheorem{proposition}[theorem]{Proposition}
\newtheorem{corollary}[theorem]{Corollaire}
\newtheorem{rappel}[theorem]{}

\theoremstyle{definition}
\newtheorem{definition}[theorem]{D\'efinition}
\newtheorem{remark}[theorem]{Remarque}

\newtheorem{notation}[theorem]{Notation}
\newtheorem{question}[theorem]{Question}
\newtheorem{example}[theorem]{Exemple}
\newtheorem{conjecture}[theorem]{Conjecture}
\newtheorem{exercice}[theorem]{Exercice}

\newcommand{\cal}{\mathcal}

\newcommand{\lan}{\langle}
\newcommand{\ran}{\rangle}

\newcommand{\cP}{\check{\P}}

\newcommand{\bi}{\begin{itemize}}
\newcommand{\ei}{\end{itemize}}
\newcommand{\be}{\begin{enumerate}}
\newcommand{\ee}{\end{enumerate}}

\newcommand{\bpf}{\begin{proof}}
\newcommand{\epf}{\end{proof}}

\newcommand{\bt}{\begin{theorem}}
\newcommand{\et}{\end{theorem}}
\newcommand{\brap}{\begin{rappel}}
\newcommand{\erap}{\end{rappel}}
\newcommand{\bnt}{\begin{notation}}
\newcommand{\ent}{\end{notation}}
\newcommand{\bd}{\begin{definition}}
\newcommand{\ed}{\end{definition}}
\newcommand{\ble}{\begin{lemma}}
\newcommand{\ele}{\end{lemma}}
\newcommand{\bpr}{\begin{proposition}}
\newcommand{\epr}{\end{proposition}}
\newcommand{\bre}{\begin{remark}}
\newcommand{\ere}{\end{remark}}
\newcommand{\bco}{\begin{corollary}}
\newcommand{\eco}{\end{corollary}}
\newcommand{\beq}{\begin{equation}}
\newcommand{\eeq}{\end{equation}}
\newcommand{\bq}{\begin{question}}
\newcommand{\eq}{\end{question}}
\newcommand{\beqn}{\begin{eqnarray*}}
\newcommand{\eeqn}{\end{eqnarray*}}
\newcommand{\bex}{\begin{example}}
\newcommand{\eex}{\end{example}}
\newcommand{\ber}{\begin{exercice}}
\newcommand{\eer}{\end{exercice}}
\newcommand{\sct}{\section}
\newcommand{\ssct}{\subsection}
\newcommand{\sk}{\smallskip}

\newcommand{\nk}{\noindent}
\newcommand{\pl}{\partial}

\newcommand{\fr}{\frac}
\newcommand{\bcj}{\begin{conjecture}}
\newcommand{\ecj}{\end{conjecture}}

\begin{document}

\begin{abstract}
Nous montrons que, si $n\geq 3$ et $d\geq 2n$, tout $d$-tissu de codimension un
pr\`es d'un point de $\C^n$, qui poss\`ede $(2d-3n+1)$ relations ab\'eliennes dont les $1$-jets 
sont lin\'eairement ind\'ependants, est isomorphe \`a un tissu 
alg\'ebrique. C'est en particulier le cas des tissus 
de rang maximal avec $n\geq 3$ et $d\geq 2n$.
Le cas $n=3$ est d\^u \`a Bol \cite{Bo}. Le cas g\'en\'eral
r\'esout un probl\`eme pos\'e par Chern et Griffiths \cite{CG1}--\cite{CG2}.

\sk

\sk

{\em We prove that, if $n\geq 3$ and $d\geq 2n$, a $d$-web in $\C^n$
is isomorphic to an algebraic $d$-web, if it has $(2d-3n+1)$
abelian relations, the $1$-jets of which are linearly independant.
The case $n=3$ is a theorem of Bol \cite{Bo}.
The general case solves a problem which was first considered by 
Chern and Giffith \cite{CG1}--\cite{CG2}.}

\end{abstract}

\title[Sur l'alg\'ebrisation des tissus]{Sur l'alg\'ebrisation des tissus, \\
le th\'eor\`eme de Bol en toute dimension $>$\! 2}

\author{Jean-Marie Tr\'epreau}

\address{UMR 7586,  175 rue du Chevaleret, 75013 Paris ; trepreau@math.jussieu.fr}

\maketitle

\include{ch-1}

\include{ch-2}

\include{ch-3}

\include{ch-4}

\include{ch-5}

\include{ch-6}

\include{bibli}
\end{document}

%% file: ch-1.tex
\sct{Introduction}

\ssct{Tissus ; tissus de rang maximal}

Pour fixer les id\'ees et les notations, nous nous pla\c{c}ons d'embl\'ee 
dans la cat\'egorie analytique complexe. Notre \'etude est locale
au voisinage de $0$ dans $\C^n$, avec $n\geq 2$. On notera $\C^n_0$ un tel 
voisinage, qu'on pourra r\'eduire autant qu'on veut. On note $(x_0,\ldots,x_{n-1})$
les coordonn\'ees d'un point de $\C^n$. Cette convention 
inhabituelle sera commode dans certains calculs.

\sk
On se donne un entier $d\geq 1$ et un $d$-tissu $\cal{T}$\footnote{On ne consid\`ere ici 
que des tissus de codimension un. Nous esp\'erons aborder le cas 
des tissus de codimension sup\'erieure dans un prochain travail.}
au voisinage de $0$ dans $\C^n$, {\em i.e.} une famille de $d$ feuilletages
de codimension un, en position g\'en\'erale. On note 
$$
u_\alpha(x), \qquad \alpha=1,\ldots,d,
$$
des fonctions de d\'efinition des feuilletages : les feuilles du $\alpha\,$-i\`eme feuilletage 
sont les hypersurfaces de niveau $\{u_\alpha(x) = {\rm cste}\}$ de la fonction 
$u_\alpha(x)$. L'hypoth\`ese de position g\'en\'erale signifie que toute 
famille d'au plus $n$ \'el\'ements extraite de la famille de
diff\'erentielles $du_1(0), \ldots, du_d(0)$ est libre.

\sk
Une relation ab\'elienne du tissu $\cal{T}$ est un $d$-uplet 
$z(x)=(z_1(x),\ldots,z_d(x))$ de fonctions analytiques
pr\`es de $0$ tel que
\beq
\label{c1-1}
\sum_{\alpha=1}^d z_\alpha(x) \, du_\alpha(x) = 0 
\eeq
et que, pour tout $\alpha$, la fonction $z_\alpha(x)$ soit constante
le long des feuilles du $\alpha\,$-i\`eme feuilletage, autrement  dit :
\beq
\label{c1-2}
dz_\alpha(x) \wedge du_\alpha(x) = 0, \qquad \alpha = 1,\ldots,d.
\eeq
Ces relations forment un espace vectoriel dont la dimension 
est appel\'ee {\em le rang du tissu}. Bol dans le cas $n=2$, puis Chern dans le cas g\'en\'eral,
ont montr\'e que le rang d'un tissu $\cal{T}$ est au plus \'egal au nombre 
entier 
\beq
\label{c1-3}
\pi(n,d) = \sum_{q=1}^{+\infty} (d - q(n-1) - 1)^+.
\eeq
(Si $m\in \Z$, on note $\,(m)^+=\max \,(0,m)$.)  Nous rappelons une d\'emonstration 
de cette in\'egalit\'e dans le \S 2.2. {\em Un tissu de rang $\pi(n,d)$ est dit de rang maximal.} 

\ssct{Le probl\`eme de l'alg\'ebrisation}

La  g\'eom\'etrie alg\'ebrique donne des exemples de tissus, 
en particulier de tissus de rang maximal.
C'est, avec le probl\`eme inverse, l'un des th\`emes majeurs 
de l'\'ecole de Blaschke entre 1927 et 1938, voir \cite{BB}.
Le probl\`eme inverse est de savoir si tous les tissus de rang maximal
viennent de la g\'eom\'etrie alg\'ebrique, \`a  diff\'eomorphisme local
pr\`es.

\sk
Soit $C$ une courbe alg\'ebrique de degr\'e $d$ dans 
l'espace projectif complexe $\P^n$, non-d\'eg\'en\'er\'ee, 
c'est-\`a-dire qui n'est pas contenue dans un hyperplan de $\P^n$.
Soit $x_0$ un point de l'espace $\cP^n$ des hyperplans de $\P^n$.
On suppose que l'hyperplan $x_0$ coupe la courbe $C$
en $d$ points distincts $p_\alpha(x_0)$, ce qui est v\'erifi\'e
pour $x_0$ g\'en\'erique. Tout hyperplan $x$ 
voisin de $x_0$ coupe aussi la courbe $C$  en $d$ points distincts
$p_\alpha(x)$. Par dualit\'e projective, les points 
$p_\alpha(x)$ d\'efinissent $d\,$ hyperplans de $\cP^n$
passant par $x$. Parce que la courbe est non-d\'eg\'en\'er\'ee, ces hyperplans sont les feuilles 
d'un $d$-tissu au voisinage de $x_0$. 

\sk
{\em Les tissus obtenus de cette mani\`ere sont dits alg\'ebriques. 
Un tissu est alg\'ebrisable s'il est
isomorphe \`a un tissu alg\'ebrique par un diff\'eomorphime
analytique local.} 

\sk
\`A ce point, le th\'eor\`eme d'addition d'Abel est fondamental.
Avec les notations pr\'ec\'edentes, si $\omega$ est une $1$-forme  holomorphe
sur $C$ et $p_0$ un point donn\'e de $C$, la somme 
$$
\sum_{\alpha=1}^d \int_{p_0}^{p_\alpha(x)} \, \omega
$$
ne d\'epend pas de l'hyperplan $x$ voisin de $x_0$. On peut en fait 
d\'emontrer que l'espace 
des $1$-formes holomorphes sur la courbe est isomorphe 
\`a l'espace des relations ab\'eliennes du tissu 
associ\'e. En particulier, le rang de ce tissu est \'egal \`a la dimension de l'espace 
des $1$-formes holomorphes sur $C$, autrement dit au genre arithm\'etique de $C$.

\sk
Le nombre $\pi(n,d)$ d\'efini par (\ref{c1-3}) est aussi 
la borne de Castelnuovo pour le genre arithm\'etique d'une 
courbe non-d\'eg\'en\'er\'ee de degr\'e $d$ dans $\P^n$. 
Une courbe non-d\'eg\'en\'er\'ee de genre $\pi(n,d)$
est dite extr\'emale. On peut ainsi associer un tissu
de rang maximal \`a toute courbe extr\'emale. 

\sk
Il y a une importante diff\'erence, qui r\'eappara\^{\i}tra,
entre le cas $n=2$ et le cas $n\geq 3$. Le 
genre arithm\'etique d'une courbe plane de degr\'e $d$ est 
toujours \'egal \`a $\pi(2,d)=(d-1)(d-2)/2$ ;
le tissu associ\'e est donc de rang maximal.
Si $n\geq 3$, la situation est toute diff\'erente. 
Les courbes extr\'emales sont l'exception.

\ssct{\'Enonc\'e du r\'esultat principal}

C'est le suivant :
\bt
\label{T1}
Si $n\geq 3$ et $d\geq 2n$, tout $d$-tissu de rang maximal au voisinage 
d'un point de $\C^n$ est alg\'ebrisable.
\et
On d\'emontre en fait un r\'esultat plus fort, voir le Th\'eor\`eme \ref{T2}.
D'autre part et modulo des variations mineures, la d\'emonstration s'applique aussi 
aux tissus r\'eels de classe $C^\infty$ au voisinage d'un point de $\R^n$.

\sk
Pour $n=3$, l'\'enonc\'e pr\'ec\'edent est un th\'eor\`eme de Bol publi\'e en 1933, 
voir \cite{Bo}. Le cas $d=2n$ est classique pour tout $n$ ; on en reparlera. 
En 1978, Chern et Griffiths ont publi\'e une d\'emonstration du 
r\'esultat g\'en\'eral, mais elle comportait une erreur qu'ils n'ont pas 
pu corriger, voir \cite{CG1}, \cite{CG2}. 
\`A notre connaissance, notre r\'esultat 
est donc nouveau quel que soit $n\geq 4$ et $d\geq (2n+1)$. 

\sk
Pour moduler cette affirmation, 
il faut dire que la nouveaut\'e,
dans la d\'emonstration, se r\'eduit \`a peu de choses. 
La strat\'egie g\'en\'erale de la d\'emonstration, qu'on 
appellera la \og m\'ethode standard \fg, est celle-l\`a m\^eme que 
Bol a introduite pour $n=3$ et que Chern et Griffiths ont ensuite 
adapt\'ee \`a la dimension quelconque. Dans le cadre de cette 
m\'ethode, le point crucial consiste \`a montrer qu'une certaine 
fonction $\phi$ de $(n+1)$  variables, construite \`a partir des relations
ab\'eliennes d'un tissu
de rang maximal et \`a valeur dans un certain espace $\P^m$,
{\em prend ses valeurs sur une surface de $\P^m$}.
La d\'emonstration de Bol est indirecte et repose 
sur une analogie entre les 
\'equations que v\'erifient les fonctions de d\'efinition
d'un tissu de rang maximal et la \og g\'eom\'etrie de Weyl \fg.
Chern et Griffiths poursuivent cette id\'ee et cherchent 
\`a montrer que ces \'equations d\'efinissent une \og g\'eom\'etrie des chemins \fg.
Notre apport consiste en la remarque qu'il s'agit,
dans les deux cas, de calculs \`a l'ordre deux et que donc, si la 
m\'ethode standard peut mener au r\'esultat escompt\'e, ce qu'on sait \^etre vrai 
au moins si $n=3$, un calcul direct doit permettre de montrer 
que l'application $\phi$ est de rang $2$ et
d'\'echapper \`a des consid\'erations
g\'eom\'etriques plus subtiles. C'est en effet le cas.
Il n'est pas exclu, mais nous ne l'avons pas v\'erifi\'e, que 
des calculs analogues permettent de terminer la 
d\'emonstration inachev\'ee de Chern et Griffiths en en conservant 
la ligne  g\'en\'erale. Quoi qu'il en soit, la d\'emonstration directe
est plus courte et plus \'el\'ementaire.

\ssct{Le th\'eor\`eme d'Abel inverse}

On appelle tissu lin\'eaire un tissu dont les feuilletages sont des 
feuilletages en (morceaux d') hyperplans. {\em Un tissu est lin\'earisable
s'il est localement diff\'eomorphe \`a un tissu lin\'eaire.}
Par construction, les tissus alg\'ebriques sont lin\'eaires.
Les tissus alg\'ebrisables sont donc lin\'earisables.

\sk
Le \og th\'eor\`eme d'Abel inverse \fg \, est un r\'esultat fondamental 
de la th\'eorie. Il s'\'enonce ainsi : tout tissu lin\'eaire qui poss\`ede 
{\em une} relation ab\'elienne $z(x)=(z_1(x),\ldots,z_d(x))$, 
dont aucune composante $z_\alpha(x)$ n'est identiquement nulle,
est alg\'ebrisable. 
Le cas $n=2$ est d\^u \`a Blaschke et Howe
et le cas g\'en\'eral fait l'objet de la premi\`ere 
partie de l'article \cite{Bo} de Bol ; voir aussi Griffiths \cite{G1}
pour des versions encore plus g\'en\'erales de ce th\'eor\`eme,
en particulier pour des tissus de codimension plus grande.
{\em Pour d\'emontrer le Th\'eor\`eme \ref{T1}, il suffit donc de montrer que, 
sous les hypoth\`eses de l'\'enonc\'e, le tissu est lin\'earisable.}

\ssct{Ce qui se passe si $\bf n=2$ ou si $\,\bf (n+1)\leq d \leq (2n-1)$}

Rappelons-le bri\`evement.
Quel que soit $n\geq 2$, un $d$-tissu est lin\'earisable si $1\leq d \leq n$ ;
il suffit de prendre ses fonctions de d\'efinition comme partie d'un syst\`eme de coordonn\'ees. 
Un $(n+1)\,$-tissu n'est pas lin\'earisable en g\'en\'eral, mais 
un $(n+1)\,$-tissu qui poss\`ede une relation ab\'elienne non triviale
(noter que $\pi(n,n+1)=1$) est lin\'earisable.
En effet, \'etant donn\'e une relation (\ref{c1-1})--(\ref{c1-2}), 
on peut, pour tout $\alpha$, \'ecrire  $z_\alpha(x)=f_\alpha(u_\alpha(x))$
et introduire une primitive $g_\alpha$ de $f_\alpha$ nulle 
en $u_\alpha(0)$. On a alors la relation :
$$
\sum_{\alpha=1}^{n+1} g_\alpha(u_\alpha(x)) = 0.
$$
On peut prendre les fonctions $g_\alpha(u_\alpha(x))$ comme fonctions de
d\'efinition des feuilletages. En changeant 
de notation, on a donc la relation $\sum_{\alpha=1}^{n+1} u_\alpha(x)  = 0$.
Il suffit maintenant de prendre $u_1(x),\ldots,u_n(x)$ 
comme nouvelles coordonn\'ees $x_0,\ldots,x_{n-1}$. Dans ce syst\`eme
de coordonn\'ees, le dernier feuilletage est d\'efini par la fonction lin\'eaire
$(x_0+\cdots+x_{n-1})$.

\sk
Du point de vue de l'alg\'ebrisation des seuls tissus de rang maximal
(on pourrait envisager des hypoth\`eses plus faibles 
que d'\^etre de rang maximal), le cas $n=2$ des tissus plans
est tr\`es particulier.
La situation est actuellement la suivante.
On a dit que tout $4$-tissu plan de rang maximal est alg\'ebrisable.
Bol a donn\'e le premier exemple d'un $d$-tissu plan de rang 
maximal non alg\'ebrisable, avec $d=5$. D'autres exemples ont \'et\'e
d\'ecouverts r\'ecemment, pour des valeurs de $d$
comprises entre $5$ et $10$ ; voir \cite{PT} et sa bibliographie.
Mis \`a part ces r\'esultats partiels, la question de l'alg\'ebrisation 
des tissus plans de rang maximal est ouverte.

\sk
On suppose maintenant $n\geq 3$ et $(n+2)\leq d\leq (2n-1)$.
Consid\'erons la famille des $d$-tissus $\cal{T}$ de la forme 
suivante. D'une part, ils contiennent le $(n+1)\,$-tissu 
$\cal{T}_0$ d\'efini par les fonctions coordonn\'ees $u_1(x)=x_0,\ldots,u_n(x)=x_{n-1}$
et la fonction 
$$
u_{n+1}(x)=x_0+ \cdots +x_{n-1}.
$$
D'autre part, leurs autres fonctions de d\'efinition 
sont de la forme 
$$
u_\alpha(x) = U_{\alpha \, 0}(x_0) + \cdots + U_{\alpha \, {n-1}}(x_{n-1}),
\qquad  \alpha = n+2,\ldots,d,
$$
en respectant seulement les conditions sur les $dU_{\alpha \, \mu}(0)$
qui assurent que les feuilletages sont en position g\'en\'erale.
Quand on les diff\'erentie, les d\'efinitions pr\'ec\'edentes
deviennent $(d-n)=\pi(n,d)$ relations ab\'eliennes ind\'ependantes : 
les  tissus de la famille sont de rang maximal par construction.

Si un diff\'eomorphisme local $\phi$ lin\'earise un tissu 
$\cal{T}$ de la famille, il transforme en particulier le $(n+1)$\,-tissu 
lin\'eaire $\cal{T}_0$ en un tissu lin\'eaire $\phi(\cal{T}_0)$.
Le tissu $\cal{T}_0$ est de rang maximal $=1$.
D'apr\`es le th\'eor\`eme d'Abel inverse, le tissu 
$\phi(\cal{T}_0)$ est associ\'e \`a une courbe alg\'ebrique 
de degr\'e $(n+1)$ dans $\P^n$. On v\'erifie d'autre part 
que les transormations qui conservent le tissu $\cal{T}_0$
sont affines. 

On obtient ainsi que les tissus lin\'earisables
de la famille ne d\'ependent que d'un nombre 
fini de param\`etres. Un tissu de la famille n'est 
donc pas lin\'earisable en g\'en\'eral.

%% file: ch-2.tex
\sct{La borne de Chern et une version pr\'ecis\'ee du Th\'eor\`eme \ref{T1}}

\ssct{Rang d'un syst\`eme de puissances de formes lin\'eaires}

On a besoin de quelques r\'esultats pr\'eliminaires sur le rang du syst\`eme des 
puissances d'ordre donn\'e d'une famille de formes lin\'eaires. 
On suppose $n\geq 2$ et $d\geq 1$. Soit 
$$
l_1(x),\ldots,l_d(x),
$$
des formes lin\'eaires en position g\'en\'erale, {\em i.e.}
telles que toute famille extraite de cardinal $\leq n$ soit libre. 
Posons :
\beq
\label{c2-1} 
r_q = \text{dim} \, (\C \, l_1(x)^q + \cdots + \C \, l_d(x)^q), \qquad q\geq 0.
\eeq
\ble
\label{C2-1}
Pour tout $n\geq 2$ et pour tout $q\geq 0$, on a les in\'egalit\'es
\beq
\label{c2-2}
r_{q+1} \geq \min(d,r_q+(n-1)), \qquad r_q\geq \min(d,q(n-1)+1).
\eeq
\ele
\bpf
Si $p\in \N$, notons $L(p)$ l'espace engendr\'e par $l_1(x)^p,\ldots,l_d(x)^p$.

On a bien s\^ur $r_0=1$ ; on a $r_1=\min(d,n)$ car les formes sont en position g\'en\'erale.
Soit maintenant $q\geq 1$. On fait les hypoth\`eses de r\'ecurrence :
$$
r_q \geq \min(d,r_{q-1}+(n-1)), \;\; r_p\geq \min(d,p(n-1)+1) \;\; \text{ si} \;\;  p=0,\ldots,q .
$$ 

\sk
On note $r=r_q$. On peut supposer que $L(q)$ est engendr\'e par $l_1(x)^q,\ldots,l_r(x)^q$.
Soit $X=\sum_{k=0}^{n-1} a_k{\pl }/\pl x_k$ un vecteur non nul tel que :
$$
X \!\cdot l_{r+\alpha} = 0, \qquad \alpha = 1,\ldots, s:=\min(d-r,n-1).
$$
L'hypoth\`ese de position g\'en\'erale permet de supposer en outre  
$$
X \!\cdot l_{\alpha}\neq 0, \qquad \alpha = 1,\ldots, r.
$$
Le vecteur $X$ induit une application lin\'eaire 
$$
X: \; L(q+1) \rightarrow L(q).
$$
Elle est surjective car 
$$
(X \!\cdot l_\alpha^{q+1})(x) = (q+1)(X \! \cdot l_\alpha) l_\alpha(x)^q
$$
avec $X \!\cdot l_\alpha\neq 0$ si $\alpha$ est compris entre $1$ et $r$. Ceci montre d\'ej\`a 
qu'on a $r_{q+1}\geq r_q$. On peut appliquer ce r\'esultat pr\'eliminaire 
aux formes $l_{r+1}(x),\ldots,l_{r+s}(x)$ :
l'espace engendr\'e par leurs puissances $(q+1)$-i\`emes 
est de dimension $\geq\min(s,q(n-1)+1)=s$.

\sk
Comme le noyau de l'application $X$ 
contient $l_{r+1}(x)^{q+1},\ldots,l_{r+s}(x)^{q+1}$, il 
est de dimension $\geq s$. On obtient 
$$
r_{q+1} \geq r + s \geq  r + \min(d-r,n-1)  \geq \min(d,r_q+(n-1))
$$
et, par hypoth\`ese de r\'ecurrence, $r_{q+1} \geq \min(d,(q+1)(n-1)+1)$.
\epf
Si $n=2$, on a ainsi $r_q\geq \min(d,q+1)$. 
Comme l'espace des polyn\^omes homog\`enes de degr\'e $q$ en deux variables 
est de dimension $(q+1)$, on en d\'eduit :
\bco
Si $n=2$, on a l'\'egalit\'e $\;r_q=\min(d,q+1)$ pour tout $q\geq 0$.
\eco
La discussion est plus subtile si $n\geq 3$. On la compl\`etera 
dans le \S 2.5 en rappelant un lemme classique de Castelnuovo.

\ssct{Relations ab\'eliennes de valuation donn\'ee}

Revenons \`a nos tissus. On suppose $n\geq 2$ et $d\geq (n+1)$.
On consid\`ere un $d$-tissu $\cal{T}$ pr\`es de $0\in\C^n$, d\'efini par des
fonctions $u_\alpha(x)$. On \'ecrit  
$$
u_\alpha(x) = u_\alpha(0) + l_\alpha(x) + O(2), \qquad \alpha=1,\ldots,d, 
$$
o\`u $l_\alpha(x)$ est une forme lin\'eaire. On note encore 
$r_q$ le rang du syst\`eme $l_1(x)^q,\ldots,l_d(x)^q$.

\sk
Soit $E(q)$ l'espace des relations ab\'eliennes
$$
z(x)=(z_1(x),\ldots,z_d(x))
$$
dont la valuation\footnote{Si $f$ est analytique au voisinage de $0$, 
on note $f=O(q)$ et on dit que $f$ est de valuation $\geq q$ si les coefficients 
des termes de degr\'e $<q$ dans le d\'eveloppement de $f$ en s\'erie 
enti\`ere sont nuls.} est $\geq q$. 
On d\'efinit ainsi une suite d\'ecroissante d'espaces ; $\,E(0)$ est l'espace 
de toutes les relations ab\'eliennes du tissu.   Soit $z(x)\in E(q)$. On \'ecrit  
$$
z_\alpha(x) = a_\alpha u_\alpha(x)^q + O(q+1), \qquad \alpha=1,\ldots,d.
$$
En prenant la partie principale de l'\'equation $\sum_{\alpha=1}^d z_\alpha(x)\, du_\alpha(x) = 0$,
on obtient :
$$
\sum_{\alpha=1}^d a_\alpha l_\alpha(x)^{q+1} = 0.
$$
On peut aussi voir la relation ci-dessus comme une \'equation lin\'eaire 
en $(a_1,\ldots,a_d)$. Il est clair que l'espace de ses solutions
est de dimension $(d-r_{q+1})$. On en d\'eduit :
\beq
\label{c2-3}
\text{dim} \, E(q) - \text{dim}\, E(q+1) \leq d-r_{q+1},
\eeq
\beq
\label{c2-4}
\text{dim} \, E(0)\leq \sum_{q=1}^{+\infty} \, (d-r_q).
\eeq
D'apr\`es le  Lemme \ref{C2-1}, on a $(d-r_q)\leq (d-q(n-1)-1)^+$. On retrouve 
ainsi l'in\'egalit\'e de Chern :  le rang $\rho$  d'un tissu est au plus \'egal au nombre 
$\pi(n,d)$ d\'efini par (\ref{c1-3}). Le lemme donne aussi :
\beq
\label{c2-5}
\rho=\pi(n,d) \; \Rightarrow \; \left( \; \forall q\geq 0, \;\; r_q = \min (d,q(n-1)+1) \; \right).
\eeq
En g\'en\'eral, l'in\'egalit\'e (\ref{c2-4}) est stricte. 

\sk
Si les fonctions de d\'efinition $u_\alpha(x)$ sont 
lin\'eaires (le tissu  est alors compos\'e de faisceaux
d'hyperplans parall\`eles), les composantes homog\`enes d'une 
relation ab\'elienne du tissu sont des relations ab\'eliennes.
On a donc :

\sk
{\em Le rang $\rho$ d'un tissu d\'efini par des formes lin\'eaires $l_1(x),\ldots,l_d(x)$
est donn\'e par la formule $\rho = \sum_{q=1}^{+\infty} \, (d-r_q)$.
Le rang d'un tissu est au plus \'egal au rang du tissu 
dont les \'el\'ements sont les faisceaux d'hyperplans parall\`eles 
aux hyperplans tangents aux feuilles du tissu initial issues d'un point donn\'e.}

\sk
Comme on a dit, la borne de Chern $\pi(n,d)$ est atteinte par les 
tissus  associ\'es aux courbes extr\'emales de $\P^n$. 
On a toutefois les exemples tr\`es simples suivants 
de $d$-tissu de rang maximal $\pi(n,d)$. Posons :
$$
l(t,x) = \sum_{\mu=0}^{n-1} t^\mu \, x_\mu.
$$
On remarque que $l(t,x)^q$ est un polyn\^ome en $t$ de degr\'e $\leq q(n-1)$.
Il en r\'esulte que, si $\theta_1,\ldots,\theta_d$ sont des nombres 
complexes deux \`a deux distincts et que l'on choisit 
$$
l_1(x) = l(\theta_1,x), \ldots, l_d(x)=l(\theta_d,x),
$$
le rang $r_q$ du syst\`eme $l_1(x)^q,\ldots,l_d(x)^q$ est 
$\leq (q(n-1)+1)$ donc $=\min(d,q(n-1)+1)$.
Le tissu d\'efini par les formes lin\'eaires 
$l_1(x),\ldots,l_d(x)$ est donc de rang $\pi(n,d)$.
Le Lemme \ref{C2-4} ci-dessous montre que la r\'eciproque est vraie : si $d\geq (2n+1)$, 
un $d$-tissu d\'efini par des formes lin\'eaires et de rang maximal 
est de la forme pr\'ec\'edente, \`a une transformation lin\'eaire de $\C^n$ pr\`es.

\ssct{Un raffinement du Th\'eor\`eme \ref{T1}}

On suppose $d\geq 2n$. Il r\'esulte de la discussion pr\'ec\'edente que la dimension 
des espaces $E(0)/E(1)$ et $E(1)/E(2)$ est major\'ee
par $(d-n)$ et $(d-2(n-1)-1)$, respectivement. La dimension de l'espace $E(0)/E(2)$ 
est donc major\'ee
par $(2d-3n+1)$. D'apr\`es (\ref{c2-5}), on a l'\'egalit\'e dans les trois cas
si le tissu est de rang maximal. 
\bd
Un $d$-tissu pr\`es d'un point de $\C^n$ est de rang maximal en valuation $\leq 1$
s'il poss\`ede un syst\`eme de $(2d-3n+1)$ relations ab\'eliennes 
dont les $1$-jets sont lin\'eairement ind\'ependants.
\ed
La condition sur les $1$-jets est \'evidemment ouverte. La discussion qui pr\'ec\`ede 
l'\'enonc\'e montre aussi qu'{\em un tissu de rang maximal en valuation $\leq 1$
a, en tout point voisin de $0$, $(d-n)$ relations ab\'eliennes dont les 
$0\,$-jets sont lin\'eairement ind\'ependants et $(d-2n+1)$ relations ab\'eliennes 
de valuation $1$, dont les $1$-jets sont lin\'eairement ind\'ependants.}
On d\'emontrera :
\bt
\label{T2}
Si $n\geq 3$ et $d\geq 2n$, un $d$-tissu de rang maximal en valuation~$\leq 1$
au voisinage d'un point de $\C^n$ est alg\'ebrisable.
\et
Le Th\'eor\`eme \ref{T1} est un cas particulier. L'id\'ee de ce 
raffinement n'est pas nouvelle. Dans le cas $n=3$, il est mentionn\'e
dans les compl\'ements au \S 35 de \cite{BB}, avec une formulation diff\'erente.

\ssct{Courbes rationnelles normales}

On se donne un entier $m\geq 2$. Les coordonn\'ees homog\`enes d'un point de $\P^m$ 
sont not\'ees $(x_0,\ldots,x_m)$. On appelle {\em points de base de $\P^m$}
les $(m+1)$ points dont toutes les coordonn\'ees sont nulles, sauf une. 
Une famille de points de $\P^m$ est en position g\'en\'erale si 
toute famille extraite de cardinal $\leq (m+1)$ est projectivement libre.

\sk
On appelle {\em courbe rationnelle normale de degr\'e $m$ dans $\P^m$}
toute courbe qui admet une param\'etrisation de la forme :
$$
x(t)=(P_0(t), \ldots, P_m(t)),
$$
o\`u $P_0,\ldots,P_m$ forment une base de l'espace des polyn\^omes de degr\'e $\leq m$
en une variable. Par la m\'ethode de Gauss, on voit qu'une courbe rationnelle
normale est projectivement \'equivalente \`a la courbe $C_m$ donn\'ee par :
$$
x(t)=(1,t, \ldots, t^{m-1},t^m).
$$
Toute famille de points d'une courbe rationnelle normale est en position 
g\'en\'erale. Il suffit de le v\'erifier pour une famille 
de $(m+1)$ points de $C_m$, autrement dit de v\'erifier que 
que si $t_0,\ldots,t_m\in \C$ sont deux \`a deux 
distincts, la relation $\sum_{j=0}^m a_j x(t_j)=0$ n'est v\'erifi\'ee 
que si tous ses coefficients sont nuls, ce qui est facile.
On a en fait la propri\'et\'e un peu plus forte suivante,
dont on laisse la d\'emonstration au lecteur :

\sk
{\em Un sous-espace $\P^d$ de dimension $d\leq (m-1)$ coupe une courbe rationnelle 
normale de $\P^m$ en au plus $(d+1)$ points, compte tenu des
multiplicit\'es.}

\sk
\nk
On a aussi :

\sk
{\em Par $(m+3)$ points de $\P^m$  en position g\'en\'erale, il passe une et 
une seule courbe rationnelle normale.}  

\sk
\nk
On peut supposer que la famille de $(m+3)$ points est constitu\'ee
des $(m+1)$ points de base et de deux autres points, repr\'esent\'es
par $x'=(x'_0,\ldots,x'_m)$ et $x''=(x''_0,\ldots,x''_m)$.

Les param\'etrages des courbes rationnelles normales qui passent par les points de base
en temps fini sont donn\'es par :
$$
x(t) = ( \prod_{j=0}^m(t-\theta_j) ) (\fr{k_0}{t-\theta_0}, \ldots, \fr{k_m}{t-\theta_m}).
$$
Un automorphisme de $\P^1$ permet encore de supposer que la courbe passe 
par les points repr\'esent\'es par $x'$ et $x''$ pour $t=\infty$ et 
$t=0$ respectivement et que $\theta_0\neq 0$ est donn\'e.
On obtient les conditions :
$$
k_j = k'x'_j, \;\;\; k_j/\theta_j=k''x''_j, \qquad j=0,\ldots,m,
$$
avec deux nouveaux param\`etres $k'$ et $k''$  non nuls.
Par homog\'en\'eit\'e, on peut supposer $k'=1$. Le premier syst\`eme 
d'\'equations d\'etermine alors $k_0,\ldots,k_m$. 
L'\'equation $k_0/\theta_0 = k''$ d\'etermine $k''$
et les autres d\'eterminent uniquement $\theta_1,\ldots,\theta_m$.

\ssct{Le lemme de Castelnuovo}

L'importance de ce lemme pour notre sujet a d'abord \'et\'e reconnue par 
Chern et Griffiths \cite{CG1}. Nous en rappelons une d\'emonstration,
inspir\'ee de \cite{CG1}. Nous continuons avec les notations 
du paragraphe pr\'ec\'edent.

\sk
Soit $\cal{O}_{m+1}(2)$ l'espace vectoriel de dimension $(m+1)(m+2)/2$
des polyn\^omes de la forme :
$$
Q(x) = \sum_{0\leq i \leq j \leq m} a_{ij} \, x_ix_j.
$$
\`A toute famille $\Gamma$ de points de $\P^m$, on associe 
le sous-espace $V(\Gamma)$ des polyn\^omes $Q\in \cal{O}_{m+1}(2)$ qui s'annulent
en tous les points de $\Gamma$.

{\em On dit que la famille $\Gamma$ impose $r$ conditions ind\'ependantes aux quadriques 
qui la contiennent si $V(\Gamma)$ est de codimension $r$ dans $\cal{O}_{m+1}(2)$.}
\ble
\label{C2-3}
Si $d\leq (2m+1)$, $\;d\,$ points en position g\'en\'erale dans $\P^m$
imposent $d$ conditions ind\'ependantes aux quadriques qui les contiennent.
\ele
\bpf
Soit $p_1,\ldots,p_d\in \P^m$ en position g\'en\'erale.
On obtient $d$ conditions sur les coefficients d'un polyn\^ome 
$Q\in V(\Gamma)$ en \'ecrivant qu'il s'annule aux points $p_1, \ldots, p_d$.
Ces conditions sont lin\'eairement ind\'ependantes. En effet,
si $i_0\in \{1,\ldots,d\}$ et compte tenu de l'hypoth\`ese de position g\'en\'erale,
il existe une r\'eunion de deux hyperplans (une quadrique !)
qui contient les $p_i$, $i\neq i_0$, mais pas $p_{i_0}$.  
\epf
Si $C$ est une courbe rationnelle normale,
param\'etr\'ee par un polyn\^ome $x(t)$ de degr\'e $m$,
$Q(x(t))$ est un polyn\^ome de degr\'e $\leq 2m$
pour tout $Q\in \cal{O}_{m+1}(2)$. On en d\'eduit que la codimension de 
$V(C)$ est $\leq (2m+1)$ donc $=(2m+1)$ d'apr\`es le lemme pr\'ec\'edent.
La r\'eciproque fait l'objet du lemme suivant :
\ble[Castelnuovo]
\label{C2-4}
Soit $d\geq (2m+3)$. Si $d\,$ points en position g\'en\'erale
dans $\P^m$ imposent $(2m+1)$ conditions ind\'ependantes 
aux quadriques qui les contiennent, ils 
appartiennent \`a une m\^eme courbe rationnelle normale de degr\'e $m$.
\ele 
\bpf
On peut supposer que la famille de points consid\'er\'ee $\Gamma$ 
contient les points de base.
D'autre part, les hypoth\`eses de l'\'enonc\'e impliquent que $\Gamma$ 
n'est pas contenu dans la r\'eunion de deux hyperplans.

Tout $Q\in V(\Gamma)$ est une combinaison lin\'eaire de  mon\^omes de la forme 
$x_i x_j$ avec $i<j$ et s'\'ecrit, de mani\`ere unique, 
$$
Q(x) = x_mL(x') + R(x'),
$$
o\`u $x'=(x_0,\ldots,x_{m-1})$, $\,L(x')$ est une forme lin\'eaire et 
$R(x')$ appartient \`a l'espace de dimension $m(m-1)/2$
engendr\'e par les mon\^omes $x_ix_j$, avec $i<j$ et $i\neq m$, $j\neq m$.
Si $R(x')$ est le polyn\^ome nul, on obtient $Q(x)=x_mL(x')$ et $Q$ doit \^etre
nul puisque $\Gamma$ n'est pas contenu dans la r\'eunion de deux hyperplans.
L'application qui associe $R(x')$ \`a $Q(x)\in V(\Gamma)$
est donc injective. Elle est aussi surjective puisque,
par hypoth\`ese, $V(\Gamma)$ est de dimension $(m+1)(m+2)/2-(2m+1)=m(m-1)/2$.

\sk
Il existe donc une base de $V(\Gamma)$ compos\'ee d'\'el\'ements de la forme :
$$
x_mL_{ij}(x') + x_ix_j, \qquad 0\leq i<j \leq m-1.
$$
On en extrait les \'el\'ements suivants :
$$
Q_j(x) = x_mL_j(x') + x_0x_j, \qquad j=1,\ldots,m-1.
$$
Pour tout $j$ et $k$, $\,j\neq k$, le  polyn\^ome
$$
L_k(x')Q_j(x)-L_j(x')Q_k(x) = x_0(L_k(x')x_j - L_j(x')x_k)
$$
s'annule sur $\Gamma$. Le polyn\^ome $L_k(x')x_j - L_j(x')x_k\in \cal{O}_{m+1}(2)$ 
s'annule aux points de $\Gamma$, sauf peut-\^etre 
aux $m$ points de base de l'hyperplan d'\'equation $\{x_0=0\}$.
Il s'annule aussi aux $(m-2)$ points de base du sous-espace 
d'\'equation $\{x_0=x_j=x_k=0\}$. Ce polyn\^ome s'annule 
donc au moins en $(d-2)\geq (2m+1)$ des points de $\Gamma$. 
Compte tenu du lemme pr\'ec\'edent et de l'hypoth\`ese,
il s'annule sur $\Gamma$ :
$$
L_k(x')x_j - L_j(x')x_k \in V(\Gamma), \qquad j,k=1,\ldots,m-1.
$$
En particulier, le coefficient de $x_j$ dans $L_k(x')$ est nul
si $j\neq k$ et $j\neq 0$. Les $L_j(x')$ sont donc de la forme suivante :
$$
L_j(x') = a_jx_j + b_jx_0, \qquad j=1,\ldots,m-1.
$$
Soit $C$ l'ensemble alg\'ebrique, qui contient $\Gamma$, 
d\'efini par les \'equations :
\begin{align*}
(a_jx_m+x_0)x_j + b_jx_mx_0                & = 0,  \qquad j=1,\ldots,m-1,                 \\
(a_k-a_j)x_jx_k + (b_kx_j - b_jx_k)x_0     & = 0, \qquad j,k=1,\ldots,m-1.
\end{align*}
Si $a_j$ ou $b_j$ \'etait nul, $a_jx_mx_j + b_jx_mx_0 + x_0x_j \in V(\Gamma)$
serait un produit de formes lin\'eaires, ce qui est impossible.
Si l'on avait $j\neq k$ et $a_j=a_k$, $\,(b_kx_j - b_jx_k)x_0\in V(\Gamma)$ serait un produit de formes lin\'eaires, 
ce qui est impossible. On a donc $a_j\neq 0$ et $b_j\neq 0$ pour tout $j$
et $a_j\neq a_k$ pour tout $j\neq k$. On en d\'eduit en particulier
que les points de $C$ tels que $x_0=0$ sont des points de base.
On param\`etre les points de $C$ tels que $x_0\neq 0$
en posant $x_0=1$, $x_m=t$ et en r\'esolvant le premier syst\`eme d'\'equations.
On obtient :
$$
x(t) = (1, -\fr{b_1 t}{1+a_1t}, \ldots, -\fr{b_{m-1} t}{1 +a_{m-1}t}, t).
$$
C'est une param\'erisation homog\`ene d'une courbe rationnelle normale de degr\'e $m$
qui passe {\em aussi} par les points de base. On conclut que cette 
courbe co\"{\i}ncide avec $C$. Le lemme est d\'emontr\'e.
\epf

%% file: ch-3.tex
\sct{Rang maximal : conditions diff\'erentielles d'ordre un, d'ordre deux}

\ssct{Introduction}

Les fonctions de d\'efinition d'un tissu de rang maximal en valuation $\leq 1$ v\'erifient 
un syst\`eme g\'en\'eral d'\'equations diff\'erentielles 
du premier et du second ordre si $n\geq 3$ et $d\geq (2n+1)$.
Ce n'est pas le cas si $n=2$. Il peut \^etre utile de pr\'eciser
que ces conditions, \`a elles seules, ne semblent pas
impliquer la lin\'earisabilit\'e. Pour d\'emontrer le th\'eor\`eme 
principal, outre ces conditions, 
nous utiliserons \`a nouveau l'existence de relations ab\'eliennes en nombre 
suffisant.  La pr\'esentation ci-dessous 
est l\'eg\`erement diff\'erente de celle de \cite{Bo}
et de celle de \cite{CG1}. 

\sk
On note $\cal{O}_n(2)$ l'espace des polyn\^omes $g$ de la forme 
\beq
\label{c3-1}
g(\xi)= \sum_{0\leq \mu\leq \nu \leq n-1} \, g_{\mu\nu} \, \xi_\mu \xi_\nu.
\eeq
Au polyn\^ome (\ref{c3-1}), on associe les op\'erateurs  \footnote{
On note $\pl_{x_\mu}$ au lieu de $\pl / \pl x_\mu$
et $\pl^2_{x_\mu x_\nu}$ au lieu de $\pl^2 / \pl x_\mu\pl x_\nu$.}
\begin{align*}
g(\nabla v)(x)  & =  \sum_{0\leq \mu\leq \nu \leq n-1} \, g_{\mu\nu} \, \pl_{x_\mu}\!v(x) \,\pl_{x_\nu}\!v(x), \\
g(\pl)v(x)      & =  \sum_{0\leq \mu\leq \nu \leq n-1} \, g_{\mu\nu} \,  \pl^2_{x_\mu x_\nu}\!v(x).
\end{align*}

\sk
On suppose $n\geq 2$ et $d\geq (2n+1)$. On consid\`ere un $d$-tissu pr\`es de $0\in \C^n$,
de rang maximal en valuation $\leq 1$, d\'efini par des fonctions $u_1(x),\ldots,u_d(x)$.

\sk
On fixe un point $x$ voisin de $0$. On \'ecrit :
$$
u_\alpha(x+y) = u_\alpha(x) + l_\alpha(y) + Q_\alpha(y) + O(3), \qquad \alpha =1,\ldots,d,
$$
o\`u $l_\alpha(y)=\sum_{\mu=0}^{n-1} y_\mu \, \pl_{x_\mu}\!u_\alpha (x)$
et $Q_\alpha (y)$ est un polyn\^ome homog\`ene de degr\'e $2$, dont les 
coefficients d\'ependent bien s\^ur de $x$. 
Soit $z=(z_1,\ldots,z_d)$ une relation ab\'elienne du tissu, avec  
$$
z_\alpha(x+y) = a_\alpha + b_\alpha l_\alpha(y) + O(2), \qquad \alpha=1,\ldots,d.
$$
Le calcul modulo $O(2)$ de la relation $\sum_{\alpha=1}^d z_\alpha \, du_\alpha=0$ donne 
\beq
\label{c3-2}
\sum_{\alpha=1}^d a_\alpha l_\alpha(y) = 0,
\qquad 
\sum_{\alpha=1}^d a_\alpha \, Q_\alpha(y) + \fr{1}{2}\sum_{\alpha=1}^d  b_\alpha l_\alpha(y)^2 =0.
\eeq
On consid\`ere maintenant (\ref{c3-2}) comme un syst\`eme d'\'equations en 
$a=(a_1,\ldots,a_d)$ et en $b=(b_1,\ldots,b_d)$. 

\ssct{Conditions du premier ordre}

On a :
\ble
\label{C3-1}
Si $n\geq 2$ et $d\geq (2n+1)$ et si le tissu $\cal{T}$ est de rang 
maximal en valuation $\leq 1$, il existe une base de $1$-formes 
$\,\omega_0(x),\ldots,\omega_{n-1}(x)$
telle que les fonctions de d\'efinition $u_1(x), \ldots,u_d(x)$
de $\cal{T}$ v\'erifient 
\beq
\label{c3-3}
du_\alpha(x) = k_\alpha(x)\sum_{\mu=0}^{n-1}  \theta_\alpha(x)^\mu \, \omega_\mu(x)
\eeq
pour des fonctions $k_\alpha(x)$ et $\theta_\alpha(x)$ convenables.
\ele
On dira d'une telle base de $1$-formes qu'elle est adapt\'ee au tissu $\cal{T}$.
\bpf
L'\'enonc\'e est invariant par diff\'eomorphisme local.
Pour obtenir l'existence d'une base de $1$-formes {\em r\'eguli\`eres}, 
il est commode de supposer que $u_1(x),\ldots,u_n(x)$ 
sont les fonctions coordonn\'ees.

Comme on l'a vu dans le \S 2.2, si le tissu est de rang maximal en valuation $\leq 1$,
il a $(d-2n+1)$ relation ab\'eliennes de valuation $1$ au point $x$, 
dont les $1$-jets sont lin\'eairement ind\'ependants. En particulier,
l'espace des solutions de la forme $(0,b)$ de (\ref{c3-2}) est de dimension $(d-2n+1)$. Comme 
$$
l_\alpha(y)^2 = \sum_{\mu,\nu=0}^{n-1} y_\mu y_\nu \,
\pl_{x_\mu}\! u_\alpha (x) \, \pl_{x_\nu}\! u_\alpha (x), \qquad \alpha=1,\ldots,d,
$$
cela revient \`a dire que la matrice de dimension $d\times n(n-1)/2$,
dont les entr\'ees sont les nombres $(\pl_{x_\mu}\! u_\alpha (x) \, \pl_{x_\nu} \!u_\alpha (x))$,
est de rang $(2n-1)$, ou encore que l'espace des polyn\^omes $g\in \cal{O}_n(2)$
tels que 
$$
g(\nabla u_\alpha)(x) = 0, \qquad \alpha=1,\ldots,d,
$$
est de codimension $(2n-1)$. 

Cette propri\'et\'e se traduit aussi de la fa\c{c}on suivante :
les $d\geq (2(n-1)+3)$ points de $\P^{n-1}$ en position g\'en\'erale rep\'esent\'es
par les vecteurs 
$$
\nabla u_\alpha(x) = (\pl_{x_0}u_\alpha(x),\ldots,\pl_{x_{n-1}}u_\alpha(x)),
$$
imposent $(2(n-1)+1)$ conditions ind\'ependantes aux quadriques 
qui les contiennent. D'apr\`es le lemme de Castelnuovo,
ils appartiennent \`a une courbe rationnelle normale de degr\'e $(n-1)$. 

La construction qu'on a d\'ecrite dans le \S 1.4 montre que la courbe 
rationnelle normale qui passe par les points de base
de $\P^{n-1}$, qui sont repr\'esent\'es par hypoth\`ese par les vecteurs
$\nabla u_1(x),\dots,\nabla u_n(x)$,  et les deux autres points repr\'esent\'es par 
$\nabla u_{n+1}(x)$ et $\nabla u_{n+2}(x)$,   admet une repr\'esentation 
param\'etrique homog\`ene de la forme 
$$
x(t) = (p_1(x,t),\ldots,p_n(x,t)),
$$
o\`u les $p_j(x,t)$ sont des polyn\^omes de degr\'e $\leq (n-1)$
en $t$ dont les coefficients d\'ependent r\'eguli\`erement de $x$.
On peut aussi \'ecrire :
$$
x(t) = \sum_{\mu=0}^{n-1} \, (a_{\mu 0}(x),\ldots,a_{\mu (n-1)}(x)) \, t^\mu.
$$
La base des $1$-formes $\omega_\mu(x) = \sum_{\lambda=0}^{n-1} a_{\mu \lambda}(x) \, dx_\lambda$
est adapt\'ee au tissu $\cal{T}$.
\epf

\ssct{Conditions du second ordre}

On continue la discussion. Le tissu $\cal{T}$ a aussi $(d-n)$ relations ab\'eliennes 
dont les $0\,$-jets au point 
$x$ sont ind\'ependants. Ceci entra\^{\i}ne en particulier que,
pour toute solution $a\in \C^d$ de la premi\`ere \'equation de (\ref{c3-2}),
il existe une solution $b\in \C^d$ de la deuxi\`eme. Donc :
$$
\sum_{\alpha=1}^d a_\alpha l_\alpha(y) = 0
\; \Rightarrow \;
\sum_{\alpha=1}^d a_\alpha \, Q_\alpha(y) \in \left( \C \, l_1(y)^2 + \cdots + \C \, l_d(y)^2 \right)
$$
pour tout $a\in \C^d$. C'est une condition sur les $Q_\alpha(y)$, non triviale si 
$l_1(y)^2,\ldots,l_d(y)^2$ n'engendrent pas l'espace des polyn\^omes 
homog\`enes de degr\'e $2$. 

Pour l'expliciter, prenons $l_1(y),\ldots,l_n(y)$ comme nouvelle base de l'espace des formes lin\'eaires.
Soit :
$$
y_\mu = \sum_{j=1}^n a_{\mu j}l_j(y), \qquad \mu = 0,\ldots, n-1,
$$
les formules de passage.
Si $\alpha\in \{1,\ldots,d\}$, on a 
$$
l_\alpha(y) = \sum_{j=1}^n  \left( \sum_{\mu=0}^{n-1} a_{\mu j} \pl_{x_\mu}\!u_\alpha(x)\right) l_j(y).
$$
On en d\'eduit :
$$
Q_\alpha(y) - \sum_{j=1}^n 
\left( \sum_{\mu=0}^{n-1} a_{\mu j} \pl_{x_\mu}\!u_\alpha(x)\right) \, Q_j(y) 
\in (\C \, l_1(y)^2+\cdots+l_d(y)^2).
$$
Soit $g\in \cal{O}_n(2)$ un polyn\^ome tel que :
$$
g(\nabla u_{\alpha'})(x)=0, \qquad \alpha'=1,\ldots,d.
$$
On remarque que :
$$
g(\pl_y)l_{\alpha'}^2(0)= \sum_{0\leq \mu\leq \nu \leq n-1} \, g_{\mu\nu} \, (\pl^2_{y_\mu y_\nu} l_{\alpha'}^2)(0)
                    = 2g(\nabla u_{\alpha'})(x) =  0,
$$
puisque $l_{\alpha'}(0)=0$, pour tout $\alpha'\in\{1,\ldots,d\}$.
On a donc $g(\pl_y)Q(0)=0$ pour tout $Q(y)\in (\C \, l_1(y)^2+\cdots+l_d(y)^2)$.
On en d\'eduit :
$$
g(\pl_y) Q_\alpha(0) = 
\sum_{j=1}^n \left( \sum_{\mu=0}^{n-1} a_{\mu j} \pl_{x_\mu}\!u_\alpha(x)\right) \, g(\pl_y) Q_j(0).
$$
D'autre part :
$$
g(\pl_y) Q_\alpha(0)=g(\pl)u_\alpha(x), 
\qquad \alpha=1,\ldots,d.
$$
En posant $m_\mu = \sum_{j=1}^n a_{\mu j}  \, g(\pl_y) Q_j(0)$, on obtient ainsi, au point {\em fix\'e} $x$ :
$$
g(\pl) u_\alpha(x) = \sum_{\mu=0}^{n-1} m_\mu \pl_{x_\mu}\! u_\alpha(x)
\qquad \alpha=1,\ldots,d.
$$
On traduit maintenant ce r\'esultat dans une base de $1$-formes 
adapt\'ee au tissu. On introduit la notation :
\beq
\label{c3-4}
\phi(x)= \sum_{\mu=0}^{n-1} (\phi)_\mu(x) \, \omega_\mu(x)  
\eeq
pour la d\'ecomposition d'une $1$-forme $\phi(x)$ dans la base $\omega_0(x),\ldots,\omega_{n-1}(x)$.
\ble
\label{C3-2}
Soit $n\geq 2$ et $d\geq (2n+1)$. Soit $\cal{T}$ un tissu de rang
maximal en valuation $\leq 1$ et $\,\omega_0(x),\ldots,\omega_{n-1}(x)$
une base de $1$-formes adapt\'ee au tissu. Avec les notations du Lemme \ref{C3-1},
on a :

\nk
Pour tout $\mu \in \{0,\ldots,n-2\}$, 
il existe des fonctions $m_{\mu 0}(x), \ldots,m_{\mu (n-1)}(x)$ telles que :
\beq
\label{c3-5}
(d(k_\alpha \theta_\alpha))_\mu(x) - (dk_\alpha)_{\mu+1}(x) 
= k_\alpha (x)\sum_{\lambda=0}^{n-1}  m_{\mu \lambda}(x) \theta_\alpha(x)^\lambda,
\qquad \alpha=1,\ldots,d.
\eeq
Si de plus $n>2$, il existe des 
fonctions $n_{\mu 0}(x), \ldots, n_{\mu n}(x)$ telles que :
\beq
\label{c3-6}
\theta_\alpha (x)(d\theta_\alpha)_\mu(x) - (d\theta_\alpha)_{\mu+1}(x)
= \sum_{\lambda=0}^n  n_{\mu \lambda}(x) \theta_\alpha(x)^\lambda,
\qquad \alpha=1,\ldots,d.
\eeq
\ele
\bre
On peut montrer que la condition (\ref{c3-6}) est une condition n\'ecessaire 
de lin\'earisabilit\'e, pour un tissu de la forme (\ref{c3-3}), m\^eme si $n=2$.
Elle n'est pas v\'erifi\'ee en g\'en\'eral par un tissu plan 
de rang maximal, mais c'est une condition n\'ecessaire (d'apr\`es ce qu'on vient de dire)
et suffisante (compte tenu de la suite de la d\'emonstration)
pour qu'un tissu plan de rang maximal soit alg\'ebrisable.
\ere
\bpf
On passe de la base canonique \`a la base adapt\'ee par des formules du type 
$$
dx_\lambda = \sum_{\lambda'=0}^{n-1} b_{\lambda \lambda'}(x) \,\omega_{\lambda'}(x),
\qquad \lambda = 0,\ldots,n-1.
$$
Pour all\'eger les \'ecritures, notons $u(x)$ l'une quelconque des fonctions $u_\alpha(x)$ et : 
$$
du(x)=  \sum_{\mu=0}^{n-1} k(x)\theta(x)^\mu \, \omega_\mu(x).
$$
On a aussi $du(x)=\sum_{\lambda=0}^{n-1} \pl_{x_\lambda}\! u(x) \, dx_\lambda$, donc :
$$
k(x)\theta(x)^\mu = \sum_{\lambda=0}^{n-1} b_{\lambda \mu}(x) \pl_{x_\lambda}\!u (x).
$$
Soit $\mu,\mu',\nu,\nu'$ quatre entiers compris entre $0$ et $(n-1)$ tels que 
$\mu + \nu=\mu'+\nu'$. De toute \'evidence,
$k(x)\theta(x)^\mu  \, k(x)\theta(x)^\nu = k(x)\theta(x)^{\mu'} \, k(x)\theta(x)^{\nu'}$.
On a donc : 
$$
\sum_{\lambda, \lambda'=0}^{n-1} 
\left( b_{\lambda \mu}(x)b_{\lambda' \nu }(x) - b_{\lambda \mu' }(x)b_{\lambda' \nu'}(x) \right)
 \,\pl_{x_\lambda}\!u (x) \, \pl_{x_{\lambda'}}\!u (x)=0.
$$
Cette relation est v\'erifi\'ee pour tout $u\in \{u_1,\ldots,u_d\}$.
Compte tenu de la discussion qui pr\'ec\`ede l'\'enonc\'e, 
il existe donc des scalaires $m_\lambda(x)$ tels que
\beq
\label{c3-7}
\sum_{\lambda, \lambda'=0}^{n-1} 
\left( b_{\lambda \mu }(x)b_{\lambda' \nu }(x) - b_{\lambda \mu'}(x)b_{\lambda' \nu'}(x) \right)
\,\pl^2_{x_\lambda x_{\lambda'}}\!u (x)
=
\sum_{\lambda=0}^{n-1} m_\lambda(x) \, \pl_{x_\lambda}\! u(x),
\eeq
pour tout $u\in \{u_1,\ldots,u_d\}$. Insistons : les $m_\lambda(x)$ d\'ependent de $\mu,\nu,\mu'$ et $\nu'$,
mais pas du choix de $u$ !

\sk
On a d'autre part :
\begin{align*}
(d(k\theta^\nu))_\mu(x) 
& = 
\sum_{\lambda=0}^{n-1} b_{\lambda \mu}(x) \, \pl_{x_\lambda}\!(k\theta^\nu)(x)  
 =
\sum_{\lambda, \lambda'=0}^{n-1} b_{\lambda \mu}(x) \, 
\pl_{x_\lambda}\!( b_{\lambda' \nu} \pl_{x_{\lambda'}}\!u  )(x)                  \\
& =
\sum_{\lambda, \lambda'=0}^{n-1} b_{\lambda \mu}(x) \,
\pl_{x_\lambda}\!b_{\lambda' \nu}(x) \, \pl_{x_{\lambda'}}\!u(x) 
+
\sum_{\lambda, \lambda'=0}^{n-1} b_{\lambda \mu}(x) \,
b_{\lambda' \nu}(x) \, \pl^2_{x_\lambda x_{\lambda'}}\!u(x).
\end{align*}
Compte tenu de (\ref{c3-7}) et du fait que $\pl_{x_\lambda}\!u(x)=k(x)\theta(x)^\lambda$,
on peut \'ecrire 
\beq
\label{c3-8}
d(k\theta^\nu))_\mu (x) -  (d(k\theta^{\nu'}))_{\mu'} (x)  
= \sum_{\lambda = 0}^{n-1}  M_\lambda (x) k(x)\theta(x)^\lambda,
\eeq
avec 
$$
M_\lambda(x)=  m_\lambda(x) + 
\sum_{\lambda'=0}^{n-1} 
( b_{\lambda' \mu}(x) \, \pl_{x_{\lambda'}}\!b_{\lambda \nu}(x) 
- b_{\lambda' \mu'}(x) \, \pl_{x_{\lambda'}}\!b_{\lambda \nu'}(x)).
$$

\sk
On choisit d'abord $\nu=1$, $\nu'=0$ et $\mu'=\mu+1$ dans (\ref{c3-8}),
ce qui donne :
$$
(d(k \theta))_\mu(x) - (dk)_{\mu+1}(x) 
= \sum_{\lambda=0}^{n-1}  m_{\mu \lambda}(x)k(x)\theta(x)^\lambda.
$$
C'est la premi\`ere partie de l'\'enonc\'e. Supposons $n\geq 3$. On peut 
alors choisir $\nu=2$, $\nu'=1$ et $\mu'=\mu+1$ dans (\ref{c3-8}), ce qui donne :
$$
(d(k\theta^2))_\mu (x) -  (d(k\theta))_{\mu+1} (x)  
= \sum_{\lambda = 0}^{n-1}  p_{\mu \lambda}(x) k(x)\theta(x)^\lambda.
$$
On a d'autre part les identit\'es :
\begin{align*}
(d(k\theta^2))_\mu(x)    & = 
\theta(x) (d(k\theta))_\mu(x) + k(x)\theta(x) (d\theta)_\mu(x),\\
(d(k\theta))_{\mu+1}(x)  & = 
\theta(x) (dk)_{\mu+1}(x) + k(x)(d\theta)_{\mu+1}(x).
\end{align*}
Par soustraction membre \`a membre et compte tenu des formules obtenues 
juste auparavant, on obtient :
$$
\sum_{\lambda = 0}^{n-1}  p_{\mu \lambda}(x) \theta(x)^\lambda
=
\theta(x)\sum_{\lambda=0}^{n-1}  m_{\mu \lambda}(x) \theta(x)^\lambda 
+
\left( \theta(x) (d\theta)_\mu(x) - (d\theta)_{\mu+1}(x) \right).
$$
C'est la deuxi\`eme partie de l'\'enonc\'e.
\epf

%% file: ch-4.tex
\sct{D\'ebut de la d\'emonstration : la m\'ethode standard}

\ssct{Introduction}

La fin de l'article est consacr\'ee \`a la d\'emonstration 
du Th\'eor\`eme \ref{T2}.
Nous suivons le plan g\'en\'eral de Bol \cite{Bo},
dont nous allons rappeler les grandes 
lignes. On pourra consulter \cite{CG1} pour une pr\'esentation 
plus d\'etaill\'ee de cette m\'ethode standard et de son arri\`ere-plan g\'eom\'etrique.

\sk
La m\'ethode repose sur des id\'ees ant\'erieures 
de Blaschke. \'Etant donn\'e un $d$-tissu de rang 
maximal $\pi(n,d):=(m+1)$, on en choisit une base de relations ab\'eliennes 
et on lui associe la famille de $d\,$  applications $p_\alpha: \C_0^n \rightarrow \P^m$,
d\'efinie par les colonnes de la matrice dont les 
lignes sont les $(m+1)$ \'el\'ements de la base. Comme $p_\alpha(x)$ ne 
d\'epend que de $u_\alpha(x)$, 
les points $p_1(x),\ldots,p_d(x)$ d\'ecrivent des courbes dans $\P^m$.
{\em Ces courbes sont des invariants projectifs du tissu} : \`a homographie pr\`es, 
elles ne d\'ependent du choix, ni des coordonn\'ees, ni des fonctions 
de d\'efinition du tissu, ni de la base de relations ab\'eliennes.
Dans \cite{B1}, Blaschke utilise cette id\'ee pour prouver
qu'un $4$-tissu plan  de rang maximal est alg\'ebrisable.
En fait, ce r\'esultat est \'equivalent \`a un r\'esultat
ant\'erieur de Lie sur les surfaces de double translation,
et Blaschke dit prendre pour mod\`ede la d\'emonstration 
qu'a donn\'ee Poincar\'e du th\'eor\`eme de Lie dans \cite{P}.
Cette premi\`ere id\'ee suffit pour d\'emontrer que les $(2n)$-tissus
de rang maximal sont alg\'ebrisables.

\sk
Si $d\geq (2n+1)$, la deuxi\`eme id\'ee essentielle consiste \`a montrer que, 
pour $x\in \C^n_0$ fix\'e,
les points $p_1(x),\ldots,p_d(x)$ appartiennent \`a une courbe rationnelle
normale $C(x)$ d'un sous-espace de $\P^m$, une propri\'et\'e qui est 
v\'erifi\'ee par les tissus qui proviennent d'une courbe alg\'ebrique
extr\'emale, et \`a \'etudier cette famille de courbes.
Cette id\'ee appara\^{\i}t d\'ej\`a dans \cite{B2}, dans une situation 
plus simple.
Il faut alors, pour mener l'analyse \`a son but, 
montrer que les courbes $C(x)$ d\'ecrivent 
une surface alg\'ebrique de $\P^m$. C'est \`a ce point 
que nous nous \'ecartons de \cite{Bo}, \cite{BB} et \cite{CG1}.
Nous n'utilisons pas l'analogie remarqu\'ee par 
Bol entre les conditons dif\'erentielles du \S 3
et l'\'equation des g\'eod\'esiques dans certaines g\'eom\'etries 
semi-riemanniennes. Nous faisons une d\'emonstration directe.

\ssct{L'application de Poincar\'e}

On suppose $n\geq 2$ et $d\geq 2n$. On consid\`ere un $d$-tissu $\cal{T}$ en $0\in \C^n$, 
de rang maximal en valuation $\leq 1$, d\'efini par des fonctions $u_1(x),\ldots,u_d(x)$.
On note :
\beq
\label{c4-1}
l: = 2d - 3n + 1, \;\;\; m := 2d - 3n.
\eeq
L'hypoth\`ese nous permet de choisir $l$ relations ab\'eliennes  
$$
z_i(x)=(z_{i\,1}(x),\ldots,z_{i\,d}(x)), \qquad i=1,\ldots l,
$$
dont les $1$-jets sont lin\'eairement ind\'ependants.
Avec un peu d'abus, on note $z'_{i,\alpha}(x)$ les fonctions d\'efinies par 
$$
dz_{i\, \alpha}(x) = z'_{i\,\alpha}(x) \, du_\alpha(x).
$$
La m\'ethode de Poincar\'e consiste d'abord \`a introduitre les vecteurs 
\beq
\label{c4-2}
Z_\alpha(x) = (z_{1\, \alpha}(x), \ldots, z_{l\, \alpha}(x)) \in \C^l,
\qquad \alpha=1,\ldots,d.
\eeq
Par d\'efinition d'une relation ab\'elienne, on a :
\beq
\label{c4-3}
\sum_{\alpha=1}^d Z_\alpha(x) \, du_\alpha(x) = 0
\eeq
et : 
$$
dZ_\alpha (x) = Z'_\alpha(x) \, du_\alpha(x), \;\; \alpha=1,\ldots,d,
$$
o\`u :
$$
Z'_\alpha(x) = (z'_{1\, \alpha}(x), \ldots, z'_{l\, \alpha}(x)).
$$
Introduisons les matrices :
$$
M(x) :=
\begin{pmatrix}
z_{1 \,1}(x)  &  \hdots  & z_{1\,d}(x)  \\
\hdots        &  \hdots  & \hdots        \\
\hdots        &  \hdots  & \hdots        \\
z_{l \,1}(x)  &  \hdots  & z_{l\,d}(x)
\end{pmatrix},
\;\; 
M'(x) :=
\begin{pmatrix}
z'_{1 \,1}(x)  &  \hdots  & z'_{1\,d}(x)  \\
\hdots         &  \hdots  & \hdots        \\
\hdots         &  \hdots  & \hdots        \\
z'_{l \,1}(x)  &  \hdots  & z'_{l\,d}(x)
\end{pmatrix}.
$$
Les $l$ lignes de $M(x)$ rep\'esentent les relations ab\'eliennes de la base qu'on a choisie
et ses $d$ colonnes repr\'esentent les vecteurs $Z_1(x),\ldots,Z_d(x)$.
\ble
\label{C4-1}
Pour tout $x$ voisin de $0$, les propri\'et\'es suivantes sont v\'erifi\'ees :

1)  les  vecteurs $Z_1(x),\ldots,Z_d(x)$ engendrent un sous-espace 
de dimension $(d-n)$  et sont en position g\'en\'erale dans ce sous-espace ;

2) toute famille \`a $(d-2n+1)$ \'el\'ements extraite
du syst\`eme $Z'_1(x),\ldots,Z'_d(x),$ engendre, ensemble 
avec les vecteurs $Z_1(x),\ldots,Z_d(x)$,  l'espace $\C^l$.
\ele 
\bpf
On peut supposer $x=0$ et un automorphisme de $\C^l$
permet de choisir la base de relations ab\'eliennes. On en construit une 
avec $(d-n)$ relations dont les $0\,$-jets sont lin\'eairement
ind\'ependants
d'une part, $\,(d-2n+1)$ relations  de valuation $1$ \`a l'origine et 
dont les $1$-jets sont lin\'eairement ind\'ependants de l'autre. Notons encore 
$$
l_\alpha(x) = \sum_{\mu=0}^{n-1} x_\mu \, \pl_{x_\mu}\! u_\alpha(0).
$$
Avec le choix qu'on a fait, les matrices $M(0)$ et $M'(0)$ sont respectivement de la forme :
$$
M(0) =
\begin{pmatrix}
A \\
O
\end{pmatrix},
\qquad
M'(0) =
\begin{pmatrix}
\star \\
B
\end{pmatrix},
$$
o\`u $O$ est une matrice nulle, $A$  une matrice de dimension 
$(d-n)\times d$ et $B$ une matrice de dimension $(d-2n+1)\times d$.
Les lignes $(a_{i\,1} \cdots a_{i \,d})$ de $A$ et 
les lignes $(b_{i\,1} \cdots b_{i \,d})$ de $B$  repr\'esentent 
respectivement des bases
de solutions des \'equations 
$$
\sum_{\alpha = 1}^d a_\alpha l_\alpha(x) = 0, 
\qquad 
\sum_{\alpha = 1}^d b_\alpha l_\alpha(x)^2 = 0.
$$
La premi\`ere assertion du lemme peut s'\'enoncer ainsi :
toutes les matrices de dimension $(d-n)\times(d-n)$ 
extraites de $A$ sont inversibles.
Si ce n'\'etait pas le cas, il existerait une combinaison lin\'eaire non triviale 
des lignes de la matrice $A$ avec au plus $n$  coefficients non nuls,
c'est-\`a-dire une relation de d\'ependance entre 
$n$ parmi les $d$ formes $l_\alpha(x)$. C'est impossible,
puisque les $l_\alpha(x)$ sont en position g\'en\'erale.

\sk
La deuxi\`eme assertion du lemme peut s'\'enoncer ainsi :
toutes les matrices de dimension $(d-2n+1)\times(d-2n+1)$ 
extraites de $B$ sont inversibles. Si ce n'\'etait pas le cas, 
il existerait une combinaison lin\'eaire non triviale 
des lignes de la matrice $B$ avec au plus $(2n-1)$ coefficients non nuls,
c'est-\`a-dire une relation de d\'ependance entre $(2n-1)$
parmi les carr\'es $l_\alpha(x)^2$. Le Lemme \ref{C2-1} montre
que c'est impossible.
\epf
On note $\C^{d-n}(x)$ le sous-espace de $\C^l$ engendr\'e par $Z_1(x),\ldots,Z_d(x)$. 
Soit 
$$
\pi : \C^l\backslash\{0\} \rightarrow \P^m
$$
la projection canonique. (Rappelons que $l=m+1$.)
Le Lemme pr\'ec\'edent montre que les applications $Z_\alpha: \C_0^n \rightarrow \C^l(x)$
induisent des applications 
\beq
\label{c4-4}
p_\alpha: \, \C_0^n \rightarrow \P^m, \qquad \alpha=1,\ldots,d,
\eeq
et que les points $p_1(x),\ldots,p_d(x)$ engendrent un sous-espace 
$\P^{d-n-1}(x)$ de dimension $(d-n-1)$. L'application induite 
\beq
\label{c4-5}
P: \, \C_0^n \rightarrow \G(d-n-1,m), \qquad P(x)=\P^{d-n-1}(x),
\eeq
\`a valeur dans la grassmannienne des $(d-n-1)$-plans de $\P^m$,
est \og l'application de Poincar\'e \fg\, associ\'ee au tissu.
\ble
\label{C4-2}
L'application de Poincar\'e (\ref{c4-5}) est une immersion.
Pour tout $x,x'\in \C_0^n$ tels que $x\neq x'$,
\beq
\label{c4-6}
\P^{n-2}(x,x'):= \P^{d-n-1}(x)\cap \P^{d-n-1}(x')
\eeq
est un sous-espace de dimension $(n-2)$ de l'espace $\P^m$.
\ele
\bpf
Soit $t\mapsto x(t)$ un arc analytique avec $x(0)=x$, $x'(0)\neq 0$. 
On a :
$$
d_t(Z_\alpha \circ x)(0) = \lan du_\alpha(x), x'(0) \ran Z'_\alpha(x), \qquad \alpha=1,\ldots,d.
$$
Le coefficient $\lan du_\alpha(x), x'(0) \ran$ est nul pour au plus $(n-1)$ 
valeurs de $\alpha$ (position g\'en\'erale).
Le Lemme \ref{C4-1} montre que les vecteurs $d_t(Z_\alpha \circ x)(0)$
et $Z_\alpha(x)$, $\alpha=1,\ldots,d$, engendrent ensemble $\C^l$.
En particulier, l'application $t\mapsto \P^{d-n-1}(x(t))$ n'est pas
stationnaire en $t=0$ : c'est la premi\`ere partie de l'\'enonc\'e.

\sk
Pour la deuxi\`eme partie, on note que, puisque l'espace engendr\'e 
par les vecteurs $(Z_\alpha\circ x)(0)$ et celui engendr\'e 
par les vecteurs $d_t(Z_\alpha \circ x)(0)$ sont transverses,
les espaces $\C^{d-n}(x(t))$ et $\C^{d-n}(x(t'))$ sont transverses 
pour tout $t,t'$ petits avec $t\neq t'$. Leur intersection est de dimension 
$2(d-n) - (2d-3n+1) = (n - 1)$. C'est la deuxi\`eme partie de l'\'enonc\'e.
\epf 

\ssct{$\bf d=2n$ ; le cas de Lie, Poincar\'e, Blaschke}

Rappelons que le fait qu'un $4$-tissu plan de rang maximal est 
alg\'ebrisable est une autre formulation d'un th\'eor\`eme difficile de 
Lie. La d\'emonstration qui suit, reprise de \cite{B1}, montre la puissance
de la m\'ethode standard, sous la forme initiale due  \`a Poincar\'e,
en liaison avec le th\'eor\`eme d'Abel inverse.

\sk
Si $d=2n$, alors $m=n$ et la grassmannienne des $(d-n-1)$-plans de $\P^n$ 
est l'espace projectif $\cP^n$ des hyperplans de  $\P^n$.
L'application de Poincar\'e 
$$
\C_0^n \rightarrow \cP^n, \qquad x\mapsto \P^{n-1}(x),
$$
est une immersion, donc un diff\'eomorphisme local.
On transporte le tissu $\cal{T}$ dans $\cP^n$ gr\^ace \`a ce diff\'eomorphisme.
Dans un syst\`eme affine de coordonn\'ees de $\cP^n$ au voisinage de $\P^{n-1}(0)$,
le tissu est lin\'eaire.
En effet, si $\alpha$ est compris entre $1$ et $d$, 
la feuille du $\alpha\,$-i\`eme feuilletage qui passe par un point 
$x_0$ est d\'efinie par l'\'equation $\{p_\alpha(x)=p_\alpha(x_0)\}$. 
Son image sous l'action 
du diff\'eomorphisme qu'on vient d'introduire est localement 
l'ensemble des hyperplans de $\P^n$ qui passent par le point 
$p_\alpha(x_0)$. Par dualit\'e projective, c'est un hyperplan.

\ssct{Les courbes rationnelles normales $\bf C(x)$}

On suppose maintenant $n\geq 2$ et $d\geq (2n+1)$. Suivant le Lemme \ref{C3-1},
on introduit une base adapt\'ee de $1$-formes $\omega_0(x),\ldots,\omega_{n-1}(x)$.
On \'ecrit : 
\beq
\label{c4-7}
du_\alpha(x) = k_\alpha(x)\sum_{\mu = 0}^{n-1}\theta_\alpha(x)^\mu \, \omega_\mu(x),
\qquad \alpha=1,\ldots,d.
\eeq
On a donc $\sum_{\alpha=1}^d Z_\alpha(x) \, du_\alpha(x) = \sum_{\mu=0}^{n-1} 
\left( \sum_{\alpha=1}^d Z_\alpha(x)k_\alpha(x)\theta_\alpha(x)^\mu \right) \,  \omega_\mu(x)$
et la relation (\ref{c4-3}) devient :
\beq
\label{c4-8} 
\sum_{\alpha=1}^d Z_\alpha(x)k_\alpha(x)\theta_\alpha(x)^\mu  = 0, \qquad \mu=0,\ldots,n-1.
\eeq

\sk
On introduit les polyn\^omes en $t\in \C$ :
$$
P(x,t) = \prod_{\beta=1}^d (t-\theta_\beta(x)), \qquad P_\alpha(x,t) = \prod_{\beta\neq \alpha}(t-\theta_\beta(x)),
$$
et la fonction $Z_\star : \C_0^n\times \C \rightarrow \C^{m+1}$ d\'efinie par :
\beq
\label{c4-9}
Z_\star(x,t) = \sum_{\alpha=1}^d \, P_\alpha(x,t) k_\alpha(x)Z_\alpha(x).
\eeq
C'est un param\'etrage homog\`ene d'une courbe rationnelle qui passe par 
les points $p_1(x),\ldots,p_d(x)$ aux temps $\theta_1(x),\ldots,\theta_d(x)$,
respectivement.  
\ble
\label{C4-3}
Les points $p_1(x),\ldots,p_d(x)$ appartiennent
\`a une courbe rationnelle normale de degr\'e $(d-n-1)$
dans l'espace $\P^{d-n-1}(x)$ qu'ils engendrent.
\ele
\bpf
On n'\'ecrit pas la variable $x$, qui est fix\'ee. 
On montre d'abord que le polyn\^ome $Z_\star(t)$,
qui est par d\'efinition de degr\'e $\leq (d-1)$, est en fait de degr\'e $\leq (d-n-1)$.

D\'efinissons les nombres $\sigma_k$ et $\sigma_k(\alpha)$ par :
$$
\prod_{\beta=1}^d (t-\theta_\beta)  = \sum_{k=0}^d \sigma_k t^k,  
\;\;\; \text{ et }  \;\;\;
\prod_{\beta\neq \alpha} (t-\theta_\beta)  = \sum_{k=0}^{d-1} \sigma_k(\alpha) t^k,
\;\;\; \alpha=1,\ldots,d.
$$
En \'ecrivant $P(t)=(t-\theta_\alpha)P_\alpha(t)$, on obtient les identit\'es
$\sigma_{k+1} = \sigma_k(\alpha) - \theta_\alpha \sigma_{k+1}(\alpha)$,
dont on tire :
$$
\sigma_k(\alpha) = \sum_{l=0}^{d-k-1} \theta_\alpha^l \, \sigma_{k+l+1}, \qquad k=0,\ldots,d-1.
$$
On peut alors calculer :
$$
Z_\star(t)   
 = \sum_{k=0}^{d-1} \left( \sum_{\alpha=1}^d \, \sigma_k(\alpha) k_\alpha \, Z_\alpha \right) \, t^k
 = \sum_{k=0}^{d-1} 
\left( \sum_{l=0}^{d-k-1} \sigma_{k+l+1} (\sum_{\alpha=1}^d  \theta_\alpha^l k_\alpha Z_\alpha ) \right)\, t^k.
$$
Dans le membre de droite, le coefficient de $t^k$ est nul 
d'apr\`es (\ref{c4-8}) si $l\leq (n-1)$ pour tout $l\leq (d-k-1)$.
Le polyn\^ome $Z_\star(t)$ est donc bien de degr\'e $\leq (d-n-1)$.

\sk
Dans un syst\`eme de coordonn\'ees de $\C^{d-n}(x)$, on a $Z_\star(t) = (Q_1(t),\ldots,Q_{d-n}(t))$,
avec des polyn\^omes $Q_1(t), \ldots, Q_{d-n}(t)$ de degr\'es $\leq (d-n-1)$. 
Comme les points $Z_1,\ldots,Z_d$ engendrent $\C^{d-n}(x)$, il est clair 
que ces polyn\^omes forment une base de l'espace des polyn\^omes de 
degr\'e $\leq (d-n-1)$. Le lemme est d\'emontr\'e.
\epf
Pour $x\in \C_0^n$, on note $C(x)$ la courbe d\'efinie par le lemme pr\'ec\'edent.
Si $x\neq x'$, les points de l'intersection $C(x)\cap C(x')$
appartiennent \`a l'espace $\P^{n-2}(x,x')$ et les rappels 
du \S 2.4 montrent qu'on a :
\ble
\label{C4-4}
Pour tout $x,x'\in \C_0^n$ tels que $x\neq x'$, les courbes $C(x)$
et $C(x')$ se rencontrent en au plus $(n-1)$ points, compte tenu des 
multiplicit\'es. 
\ele 
On peut affirmer que $C(x)$ et $C(x')$ se coupent en exactement 
$(n-1)$ points dans le cas suivant : $u_\alpha(x)=u_\alpha(x')$
pour exactement $(n-1)$ valeurs de l'indice $\alpha$.
Les courbes ont alors les points correspondants $p_\alpha(x)=p_\alpha(x')$
en commun. {\em On montrera plus bas que deux 
courbes $C(x)$ distinctes se coupent toujours en exactement $(n-1)$
points, compte tenu des multiplicit\'es.}

%% file: ch-5.tex
\sct{Suite de la d\'emonstration : les courbes C(x) engendrent une surface}

\ssct{\'Enonc\'e du lemme principal}

C'est le suivant :
\ble
\label{C5-1}
On suppose $n\geq 3$ et $d\geq (2n+1)$. L'application 
\beq
\label{c5-1}
p: \C^n_0 \times \P^1 \rightarrow \P^m,
\eeq
induite par (\ref{c4-9}) est de rang $2$ en tout $(x,t)\in \C_0^n\times \P^1$. 
\ele
L'\'enonc\'e est vrai si $n=2$, \`a condition de supposer que la deuxi\`eme propri\'et\'e
du Lemme \ref{C3-2} est v\'erifi\'ee, ce qui dans ce cas n'est pas impliqu\'e
par les autres hypoth\`eses.

\sk
L'image de l'espace tangent \`a $\C_0^n\times \C$ en un point $(x,t)$ par l'application 
d\'eriv\'ee de $p$ est la projection, par la d\'eriv\'ee de la
projection canonique $\pi : \C^{m+1}\backslash\{0\}\rightarrow \P^m$, du sous-espace de $\C^{m+1}$
engendr\'e par les vecteurs :
$$
Z_\star(x,t), \;\; \pl_t Z_\star(x,t), \;\; \pl_{x_0}Z_\star(x,t), \ldots, \pl_{x_{n-1}}Z_\star(x,t).
$$
Il s'agit de montrer que cet espace est de dimension $3$. Introduisons la $1$-forme suivante :
\beq
\label{c5-2}
\Omega (x,t):= \sum_{\mu=0}^{n-1} t^\mu \, \omega_\mu(x), 
\qquad (x,t)\in \C_0^n\times \C.
\eeq
Le r\'esultat cherch\'e est alors une cons\'equence imm\'ediate du suivant.
\ble
\label{C5-2}
On suppose $n\geq 3$ et $d\geq (2n+1)$. Il  existe un vecteur $F(x,t)$ 
et des $1$-formes $\Gamma(x,t)$ et $\Delta(x,t)$ 
telles que\footnote{ La notation $dZ_\star(x,t)$ renvoie \`a la diff\'erentielle 
de $Z(x,t)$ en $x$. Il en ira de m\^eme dans la suite.} :
\beq
\label{c5-3}
dZ_\star(x,t) = \Omega(x,t) \,F(x,t)  + \Gamma(x,t) \, Z_\star(x,t)  + \Delta(x,t) \, \pl_t Z_\star(x,t).
\eeq
De plus les vecteurs $F(x,t)$, $Z_\star(x,t)$ et $\pl_t Z_\star(x,t)$ sont lin\'eairement 
ind\'ependants. 
\ele
Ce lemme montre que si $s\mapsto x(s)$ est un arc 
tel que $x(0)=x$, l'arc $s\mapsto p(x(s),t)$ est transverse \`a la courbe 
$C(x)$ au point $p(x,t)$ si et seulement si $x'(0)$ 
n'appartient pas au noyau de la forme lin\'eaire $\Omega(x,t)$.

\ssct{D\'ebut de la d\'emonstration du Lemme \ref{C5-2}}

On n'\'ecrit plus les variables $(x,t)$. Plut\^ot qu'avec $Z_\star$, on 
travaille avec la fonction suivante :
$$
Z = \sum_{\alpha=1}^d \fr{k_\alpha\, Z_\alpha}{t-\theta_\alpha}.
$$
Rappelons que la formule $dZ = \sum_{\mu=0}^{n-1} (dZ)_\mu\, \omega_\mu$
d\'efinit les fonctions $(dZ)_\mu$. On a :
$$
dZ = \sum_{\alpha=1}^d \, \fr{k_\alpha\, Z'_\alpha}{t-\theta_\alpha} \,du_\alpha 
   + \sum_{\alpha=1}^d \, Z_\alpha \, d\left( \fr{k_\alpha}{t-\theta_\alpha} \right)
$$
Compte tenu de la forme (\ref{c4-7}) des $du_\alpha$, on a donc : 
\beq
\label{c5-4}
(dZ)_\nu  = \sum_{\alpha=1}^d \, \fr{k_\alpha^2\theta_\alpha^\nu \, Z'_\alpha}{t-\theta_\alpha}   \, + \,
\sum_{\alpha=1}^d \, Z_\alpha \, \left( d\left( \fr{k_\alpha}{t-\theta_\alpha} \right) \right)_\nu,
\qquad \nu=0,\ldots,n-1.
\eeq
En diff\'erentiant la relation $\sum_{\alpha=1}^d Z_\alpha k_\alpha = 0$,
voir (\ref{c4-8}), on obtient la relation 
$$
\sum_{\alpha=1}^d Z'_\alpha k_\alpha \, du_\alpha + \sum_{\alpha=1}^d Z_\alpha \, dk_\alpha = 0,
$$
qu'on d\'ecompose dans la base adapt\'ee :
\beq
\label{c5-5}
\sum_{\alpha=1}^d Z'_\alpha k_\alpha^2 \theta_\alpha^\mu  = -\sum_{\alpha=1}^d Z_\alpha \, (dk_\alpha)_\mu,
\qquad \mu=0,\ldots,n-1.
\eeq
On peut s'\'etonner qu'on utilise une seule des relations (\ref{c4-8}) : 
les autres sont en fait couvertes par le Lemme \ref{C3-2}, de m\^eme d'ailleurs que 
celles qu'on obtiendrait en \'ecrivant que le second membre de 
(\ref{c4-7}) est une forme ferm\'ee.

\ssct{Un calcul}

La cl\'e du calcul est le lemme suivant :
\ble
Il existe des fonctions $f_\mu(x,t)$ et $g_\mu(x,t)$ telles que :
\beq
\label{c5-6}
t(dZ)_\mu - (dZ)_{\mu+1}= f_\mu \, Z + g_\mu \, \pl_t Z,
\qquad \mu=0,\ldots,n-2.
\eeq
\ele
\bpf
Dans ce qui suit, l'entier $\mu\in \{0,\ldots,n-2\}$ est fix\'e. Notons :
$$
I_\mu : = t(dZ)_\mu - (dZ)_{\mu+1}.
$$
Les formules (\ref{c5-4}) donnent :
$$
I_\mu =  \sum_{\alpha=1}^d \, k_\alpha^2\theta_\alpha^\mu \, Z'_\alpha   \, + \,
\sum_{\alpha=1}^d \, Z_\alpha \, 
\left( t\left( d\left( \fr{k_\alpha}{t-\theta_\alpha} \right) \right)_\mu
-
\left( d\left( \fr{k_\alpha}{t-\theta_\alpha} \right) \right)_{\mu+1} \right).
$$
Les formules (\ref{c5-5}) permettent d'\'eliminer les $Z'_\alpha$. On obtient ainsi : :
$$
I_\mu = \sum_{\alpha=1}^d Z_\alpha \, K_{\mu \alpha},
$$
o\`u (rappelons que les diff\'erentielles concernent les $x$, pas les $t$) :
\begin{align*}
K_{\mu \alpha} 
& =
- (dk_\alpha)_\mu  + 
t \left( d\left(\fr{k_\alpha}{t-\theta_\alpha}\right)\right)_\mu 
- \left( d\left(\fr{k_\alpha}{t-\theta_\alpha}\right)\right)_{\mu + 1} \\
& = 
 \left( d\left( \fr{k_\alpha\theta_\alpha}{t-\theta_\alpha}\right)\right)_\mu 
- \left( d\left(\fr{k_\alpha}{t-\theta_\alpha}\right)\right)_{\mu + 1}     \\
& = 
 \fr{ (d(k_\alpha\theta_\alpha))_\mu - (dk_\alpha)_{\mu+1}}{t-\theta_\alpha}
+
 \fr{ k_\alpha (\theta_\alpha (d\theta_\alpha)_\mu - (d\theta_\alpha)_{\mu+1}}{(t-\theta_\alpha)^2}\\
& =
L_{\mu \alpha} + M_{\mu \alpha}.
\end{align*}
C'est (enfin) le moment de rappeler les formules suivantes du Lemme \ref{C3-2} :
\begin{align*}
(d(k_\alpha \theta_\alpha))_\mu - (dk_\alpha)_{\mu+1} 
& = \sum_{\lambda=0}^{n-1}  m_{\mu \lambda} k_\alpha \theta_\alpha^\lambda, \\
\theta_\alpha (d\theta_\alpha)_\mu - (d\theta_\alpha)_{\mu+1}
& = \sum_{\lambda=0}^n n_{\mu \lambda} \theta_\alpha^\lambda.
\end{align*}
On a donc :
$$
L_{\mu \alpha} 
 =  \sum_{\lambda=0}^{n-1}  \fr{m_{\mu \lambda} k_\alpha \theta_\alpha^\lambda}{t-\theta_\alpha}
 =  \sum_{\lambda=0}^{n-1}  \fr{m_{\mu \lambda} k_\alpha t^\lambda}{t-\theta_\alpha}
 - \sum_{\lambda=1}^{n-1}\sum_{\lambda'=0}^{\lambda-1} 
m_{\mu \lambda }k_\alpha t^{\lambda-\lambda'} \theta_\alpha^{\lambda'}. 
$$
De la m\^eme mani\`ere, on \'ecrit :
\begin{align*}
M_{\mu \alpha}
& = \sum_{\lambda=0}^n \fr{n_{\mu \lambda} k_\alpha \theta_\alpha^\lambda}{(t-\theta_\alpha)^2} 
  = 
\sum_{\lambda=0}^n \fr{n_{\mu \lambda} k_\alpha t^\lambda}{(t-\theta_\alpha)^2} 
-
\sum_{\lambda=1}^n \sum_{\lambda'=0}^{\lambda-1}
\fr{n_{\mu \lambda} k_\alpha t^{\lambda-\lambda'}\theta_\alpha^{\lambda'}}{t-\theta_\alpha} \\
& =
\sum_{\lambda=0}^n \fr{n_{\mu \lambda} k_\alpha t^\lambda}{(t-\theta_\alpha)^2} 
-
\sum_{\lambda=1}^n \fr{\lambda n_{\mu \lambda} k_\alpha t^\lambda}{t-\theta_\alpha}
+
\sum_{\lambda=2}^n \sum_{\lambda'=1}^{\lambda-1}\sum_{\lambda''=0}^{\lambda'-1}
n_{\mu \lambda} k_\alpha t^{\lambda-\lambda''}\theta_\alpha^{\lambda''}.
\end{align*}
On a donc, compte tenu de (\ref{c4-8}) :
$$
\sum_{\alpha=1}^d L_{\mu \alpha}\,Z_\alpha =  \left( \sum_{\lambda=0}^{n-1} m_{\mu \lambda} t^\lambda \right) Z,
$$
$$
\sum_{\alpha=1}^d M_{\mu \alpha}\, Z_\alpha =
-   \left( \sum_{\lambda=0}^n n_{\mu \lambda} t^\lambda \right) \,  \pl_t Z
-   \left( \sum_{\lambda=1}^n \lambda n_{\mu \lambda} t^\lambda \right) Z,
$$
ce qui termine la d\'emonstration du lemme.
\epf

\ssct{Fin de la d\'emonstration du Lemme \ref{C5-2}}

On a montr\'e :
$$
(dZ)_\mu = t(dZ)_{\mu-1} - f_{\mu-1} \, Z - g_{\mu-1} \, \pl_t Z, \qquad \mu=1,\ldots,n-1.
$$
On en d\'eduit 
$$
(dZ)_\mu =  t^\mu (dZ)_0 - Z\,       \sum_{\lambda=0}^{\mu-1} t^{\lambda}f_{\mu-\lambda-1}
                         - \pl_t Z \, \sum_{\lambda=0}^{\mu-1} t^{\lambda}g_{\mu-\lambda-1},
 \qquad \mu=1,\ldots,n-1,
$$
et une d\'ecomposition de la forme :
$$
dZ(x,t) = \Omega(x,t) \,F'(x,t)  + \Gamma'(x,t) \, Z(x,t)  + \Delta'(x,t) \, \pl_t Z(x,t).
$$
Comme $Z_\star = PZ$, on en d\'eduit (\ref{c5-3}) avec 
$$
F(x,t) = P(x,t)(dZ)_0(x).
$$
Il reste \`a montrer que $F(x,t)$ n'appartient \`a l'espace engendr\'e par 
$Z(x,t)$ et $\pl_t Z(x,t)$. On a en fait le r\'esultat plus fort suivant :
$F(x,t)$ n'appartient pas \`a l'espace $\C^{d-n}(x)$
engendr\'e par $Z_1(x),\ldots,Z_d(x)$.

\sk
On a :
$$
(dZ)_0(x,t)  = \sum_{\alpha=1}^d \, \fr{k_\alpha(x)^2 \, Z'_\alpha(x)}{t-\theta_\alpha(x)} 
             + \sum_{\alpha=1}^d \, Z_\alpha(x) \, \left( d\left( \fr{k_\alpha(x)}{t-\theta_\alpha(x)} \right) \right)_0.
$$
Comme le second terme de la somme appartient \`a $\C^{d-n}(x)$, il suffit de montrer que 
$$
X(x,t)  = \sum_{\alpha=1}^d \, P_\alpha(x,t) k_\alpha(x)^2 \, Z'_\alpha(x) 
$$
n'appartient pas \`a $\C^{d-n}(x)$. 
On \'ecrit, avec les notations introduites pour d\'emontrer le Lemme \ref{C4-3}  :
\begin{align*}
X(x,t)  & = \sum_{\alpha=1}^d \sum_{k=0}^{d-1} \, \sigma_k(\alpha)t^k k_\alpha(x)^2 \, Z'_\alpha(x)         \\
        & = \sum_{\alpha=1}^d \sum_{k=0}^{d-1} \, 
            \sum_{l=0}^{d-k-1} \theta_\alpha^l \, \sigma_{k+l+1} t^k \, k_\alpha(x)^2 \, Z'_\alpha(x).
\end{align*}
En diff\'erentiant les relations (\ref{c4-8}), on obtient que, si $\mu\in\{0,\ldots,n-1\}$,
$$
\sum_{\alpha=1}^d Z'_\alpha(x)k_\alpha(x)\theta_\alpha(x)^\mu \, du_\alpha(x) \in \C^{d-n}(x).
$$
En \'ecrivant ceci dans la base adapt\'ee, on obtient la m\^eme propri\'et\'e pour les sommes :
$$
\sum_{\alpha=1}^d Z'_\alpha(x)k_\alpha(x)^2\theta_\alpha(x)^\mu \, du_\alpha(x), \qquad \mu=0,\ldots,2n-2.
$$

\sk
On en d\'eduit que $X(x,t)$ est la somme d'un \'el\'ement de $\C^{d-n}(x)$
et d'un polyn\^ome en $t$ de degr\'e $\leq d-2n$.
Le Lemme \ref{C4-1} nous dit par ailleurs que les classes modulo $\C^{d-n}(x)$
des vecteurs $Z'_1(x),\ldots,Z'_d(x)$ engendrent un espace de 
dimension $(d-2n+1)$. Ceci implique que la classe de $X(x,t)$ modulo $\C^{d-n}(x)$
ne s'annule pas. Le lemme est d\'emontr\'e.

%% file: ch-6.tex
\sct{Fin de la d\'emonstration}

\ssct{Introduction}

Arm\'es du Lemme \ref{C5-2}, nous nous retrouvons en terrain connu.
Les articles d\'ej\`a cit\'es \cite{Bo}, \cite{BB} et \cite{CG2}
fournissent plusieurs mani\`eres d'achever la lin\'earisation 
du tissu $\cal{T}$.

\sk
Dans leur livre \cite{BB}, Blaschke et Bol indiquent une voie que nous 
allons suivre. R\'esumons-la bri\`evement.

\sk
On sait maintenant que les courbes $C(x)$ sont port\'ees par une surface 
$S_0\subset \P^m$. On d\'eduira d'abord, du fait
que $S_0$ contient beaucoup de courbes rationnelles 
normales, que $S_0$ est contenu dans une surface alg\'ebrique 
$S$ de $\P^n$.
On montrera ensuite que la famille des courbes $C(x)$ est contenue 
dans une famille alg\'ebrique $\cal{C}$, de dimension $n$, de courbes  
alg\'ebriques port\'ees par $S$. 
Compte tenu des propri\'et\'es d'intersection 
des courbes $C(x)$ et par continuit\'e, on obtiendra que,
par toute famille g\'en\'erique de $n$ points de $S$, il passe une et une 
seule courbe de la famille $\cal{C}$. Nous serons alors en mesure 
d'appliquer un th\'eor\`eme classique d'Enriques et de conclure que 
le syst\`eme $\cal{C}$ est un syst\`eme lin\'eaire.
Autrement dit, il est param\'etr\'e par l'espace projectif $\P^n$,
de telle fa\c{c}on que l'ensemble 
de ses \'el\'ements qui passent par un point donn\'e 
de $S$ est param\'etr\'e par un hyperplan de $\P^n$.
On se retrouvera ainsi dans une situation o\`u l'argument 
utilis\'e dans le cas $d=2n$ pourra s'appliquer :
le diff\'eomorphisme local $\C_0^n \rightarrow \P^n$
qui envoie $x$ sur le param\`etre de la courbe $C(x)$
lin\'earise le tissu.
Compte tenu du th\'eor\`eme d'Abel inverse, le tissu $\cal{T}$
est alg\'ebrisable.
\bre
La d\'emonstration initiale de Bol \cite{Bo} n'utilise pas le 
th\'eor\`eme d'Enriques, mais des projections successives 
bien choisies, jusqu'\`a se ramener \`a la situation qu'on 
rencontre dans le cas $d=2n$. L'hypoth\`ese que le tissu 
est de rang maximal, et pas seulement de rang maximal en valuation 
$\leq 1$, est alors importante, au moins si l'on suit 
la d\'emonstration de Bol pas \`a pas. 
Il n'est pas exclu que cette m\'ethode plus 
\'el\'ementaire puisse \^etre adapt\'ee au cas g\'en\'eral.
\ere 

\ssct{Le point sur la situation}

Il est peut \^etre utile de rappeler o\`u nous en sommes dans la
d\'emonstration. On suppose toujours $d\geq (2n+1)$ et $n\geq 3$\footnote{
Le cas $n=2$ est permis \`a condition de faire l'hypoth\`ese
suppl\'ementaire que la conclusion du Lemme \ref{C3-2}
est v\'erifi\'ee. Rappelons aussi que $\C_0^n$ d\'esigne un voisinage
de $0$ dans $\C^n$,
qu'on peut restreindre autant qu'on veut.}.

\sk
On consid\`ere un $d$-tissu $\cal{T}$ au voisinage de $0\in \C^n$,
qui poss\`ede $m+1=(2d-3n+1)$ relations ab\'eliennes dont les $1$-jets sont lin\'eairement 
ind\'ependants.

\sk
On a associ\'e au tissu $\cal{T}$ une famille d'applications :
$$
p_\alpha: \C_0^n \rightarrow \P^m, \qquad \alpha = 1,\ldots,d.
$$
L'image de chacune d'elles est une courbe dans $\P^m$.
La feuille $\cal{F}_\alpha(x_\star)$ du $\alpha-$i\`eme feuilletage 
qui passe par un point $x_\star$ est donn\'ee par 
$$
\cal{F}_\alpha(x_\star) = \{x\in \C_0^n,  \;\; p_\alpha(x)=p_\alpha(x_\star)\}, 
\qquad \alpha=1,\ldots,d.
$$
On sait que, pour tout $x$, les points $p_1(x),\ldots,p_d(x)$ engendrent un 
sous-espace $\P^{d-n-1}(x)$ de dimension $(d-n-1)$ et qu'ils sont 
en position g\'en\'erale dans ce sous-espace. De plus, l'application
de Poincar\'e ainsi d\'efinie 
$$
\P^{d-n-1}: \; \C_0^n \rightarrow \G(d-n-1,n)
$$ 
est une immersion. Pour tout couple de points distincts $x,x'\in \C_0^n$,
les sous-espaces $\P^{d-n-1}(x)$ et $\P^{d-n-1}(x')$ se coupent (transversalement)
suivant un sous-espace $\P^{n-2}(x,x')$ de dimension $(n-2)$.

\sk
On sait que, pour tout $x$, les points 
$p_1(x),\ldots,p_d(x)$ appartiennent \`a une courbe rationnelle normale 
$C(x)$ de $\P^{d-n-1}(x)$, uniquement d\'etermin\'ee. Enfin, d'apr\`es le 
Lemme \ref{C5-1}, il existe une application de rang constant $2$ :
$$
p: \; \C_0^n\times \P^1 \rightarrow \P^m,
$$
telle que, pour tout $x\in \C_0^n$, $\;t\mapsto p(x,t)$ est un isomorphisme de $\P^1$ sur $C(x)$.

\sk
D'apr\`es le th\'eor\`eme du rang, l'image de $p$ est une sous-vari\'et\'e lisse de $\P^m$, non ferm\'ee.
On note cette surface $S_0$. 

\ssct{Propri\'et\'es d'intersection des courbes $C(x)$}

On a :
\ble
\label{C6-1}
Pour tout couple de points distincts $x,x'\in \C_0^n$,
les courbes $C(x)$ et $C(x')$ se coupent en exactement 
$(n-1)$ points, compte tenu des multiplicit\'es.
\ele 
\bpf
On a vu que $C(x)$ et $C(x')$ se coupent en au plus 
$(n-1)$ points.
(Rappelons que c'est une cons\'equence du fait que $\P^{d-n-1}(x)$
et $\P^{d-n-1}(x')$ sont transverses d'intersection $\P^{n-2}(x,x')$
et que toute famille de points d'une courbe rationnelle normale est 
en position g\'en\'erale dans l'espace qu'elle engendre.)
Pour $t\in \C$ voisin de $0$, soit $x(t)\in \C_0^n$
le point d\'efini par $u_1(x(t))=0, \ldots, u_{n-1}(x(t))=0$ et 
$u_n(x(t))=t$. La courbe $C(x(t))$ contient les points 
$p_1(0),\ldots,p_{n-1}(0)$, qui sont en position g\'en\'erale dans $\P^m$,
donc distincts. Ainsi, pour tout $t\neq 0$ petit, 
$C(x(t))$ coupe $C(0)$ en exactement $(n-1)$ points. 
{\em Compte tenu du fait que les courbes $C(x)$ sont contenues 
dans une surface}, le nombre de points d'intersection de deux 
courbes distinctes $C(x)$ et $C(x')$ est constant quand $x'$ varie un peu.
D'o\`u le lemme.
\epf

\ssct{La surface $S_0$ est contenue dans une surface alg\'ebrique $S$}

En effet, donnons-nous un point de $S_0$
 et un voisinage $S_1$ de ce point dans $S_0$,
contenu dans 
le domaine d'une carte affine de $\P^m$. On se ram\`ene ainsi au cas o\`u 
$0\in S_1\subset \C^m$. Il suffit de montrer que $S_1$
est contenu dans une surface alg\'ebrique de $\C^m$.

\sk
Au voisinage de $0$, $S_1$ est d\'efini par un syst\`emes d'\'equations,
soit $f_\rho(x)=0$ avec $\rho=1,\ldots,r$.
Consid\'erons l'espace $X$ des polyn\^omes $x: \C\rightarrow \C^m$
de degr\'e $\leq (d-n-1)$ et nuls en $0\in \C$.
Un \'el\'ement de $X$ est donn\'e par la 
famille convenablement ordonn\'ee de ses 
coefficients, soit $a=(a_1,\ldots,a_N)$. On note $x_a$ 
le polyn\^ome donn\'e par $a\in \C^N$. 

\sk
Si $\rho\in \{1,\ldots,r\}$, quand on  d\'eveloppe $f_\rho(x_a(t))$
en s\'erie enti\`ere de $t$, les coefficients de la 
s\'erie obtenue sont des polyn\^omes en $a$.
L'ensemble $X'$ des $x\in X$ qui envoie un voisinage de 
$0\in \C$ dans $S_1$ est donc un ensemble alg\'ebrique.
Il contient des param\'etrisations
convenables des courbes $C(x)$ qui passent par $0$.
L'image de l'application $(a,t)\mapsto x_a(t)$ est une 
surface alg\'ebrique de $\C^m$, qui contient $S_0$. 

\ssct{Le syst\`eme $\cal{C}_0$ des courbes $C(x)$ est contenu dans un syst\`eme
alg\'ebrique de dimension $n$}

On peut munir l'ensemble $\cal{R}$ des courbes rationnelles normales 
de degr\'e $(d-n-1)$ dans $\P^m$ d'une structure naturelle de vari\'et\'e analytique lisse.
Comme tous les courbes $C\in \cal{R}$ sont \'equivalentes,
il suffit de d\'efinir la structure analytique au voisinage 
d'une courbe, par exemple la courbe $C_0$ param\'et\'ee par :
$$
x(t)=(1,t,\ldots,t^{d-n-1},0,\ldots,0).
$$
On peut alors introduire une famille suffisante d'hyperplans $H_k$
transverses \`a la courbe $C_0$ et param\'etrer les courbes 
voisines de $C_0$ par leurs intersections avec ces hyperplans...
Les d\'etails sont laiss\'es au lecteur. 

\sk
Il r\'esulte du fait que l'application de Poncar\'e $\P^{d-n-1}: \C_0^n \rightarrow \G(d-n-1,m)$
est une immersion que l'application $x\mapsto C(x)$ d\'efinit une immersion :
\beq
\label{c6-1}
C: \;  \C_0^n \rightarrow \cal{R}.
\eeq
On note $\cal{C}_0$ son image.

\sk
La th\'eorie des vari\'et\'es de Chow montre que la vari\'et\'e 
$\cal{R}$ est une vari\'et\'e quasi-projective.
Autrement dit, il existe une sous-vari\'et\'e (ferm\'ee) $\cal{R}'$
d'un espace projectif, telle que $\cal{R}$ soit isomorphe 
\`a un ouvert de Zariski dense de $\cal{R}'$.
On identifie provisoirement $\cal{R}$ \`a son image 
dans $\cal{R}'$. 
Plus pr\'ecis\'ement, il existe aussi une sous-vari\'et\'e 
$\cal{W}$ de $\cal{R}'\times \P^m$, de dimension 
$$
\text{dim} \, \cal{W} = \text{dim} \, \cal{R} + 1,
$$
qu'on munit des projections canoniques  
$$
\varpi: \; \cal{W} \rightarrow \cal{R}', \qquad 
\pi: \;    \cal{W} \rightarrow \P^m,
$$
et qui a la propri\'et\'e suivante. 
Pour tout $\xi\in \cal{R}'$, son image r\'eciproque 
$\varpi^{-1}(\xi)$ est une courbe de $\cal{W}$.
De plus, si $\xi\in \cal{R}$, la projection 
$\pi$ induit un isomorphisme de cette courbe 
sur son image dans $\P^m$ et cette image est la courbe rationnelle 
normale repr\'esent\'ee par $\xi$.
\ble
\label{C6-3}
La famille $\cal{C}_0$ des courbes $C(x)$, $x\in \C_0^n$, 
est contenue dans une sous-vari\'et\'e alg\'ebrique irr\'eductible $\cal{C}$
de $\cal{R}'$, de dimension $n$.
\ele
\bpf
Suivant \cite{BB}, on consid\`ere l'intersection $\cal{C}$ de toutes 
les sous-vari\'et\'es alg\'ebriques de $\cal{R}'$ qui contiennent $\cal{C}_0$. 
C'est une vari\'et\'e irr\'eductible qui contient $\cal{C}_0$
et qui est donc de dimension $\geq n$. Il s'agit de montrer qu'elle est 
de dimension $n$. 

\sk
Si $\xi\in \cal{R}'$, on note aussi $\xi$ la courbe $\pi(\varpi^{-1}(\xi))$.
L'ensemble des $\xi\in \cal{C}$ qui sont contenus dans $S$
est une sous-vari\'et\'e de $\cal{C}$  qui contient $\cal{C}_0$.
C'est donc $\cal{C}$.

\sk
Pour tout $x\in \C_0^n$, l'ensemble des $\xi \in \cal{C}$
qui coupe $C(x)$ en $(n-1)$ points est un ouvert de Zariski
d'une sous vari\'et\'e qui contient $\cal{C}_0$, donc un ouvert 
de Zariski de $\cal{C}$.

\sk
De m\^eme, si $\xi$ est un \'el\'ement de $\cal{C}$ qui coupe toutes les 
$C(x)$ en $(n-1)$ points, l'ensemble des $\xi'\in \cal{C}$ qui coupent 
$\xi$ en $(n-1)$ points est un ouvert
de Zariski d'une sous-vari\'et\'e de $\cal{C}$
qui contient $\cal{C}_0$, donc un ouvert de Zariski de $\cal{C}$.

\sk
Finalement, si la dimension de $\cal{C}$ \'etait 
$\geq (n+1)$, deux \'el\'ements g\'en\'eriques de
$\cal{C}$ se couperaient en au moins $n$ points, ce qui 
contredit ce qu'on vient de montrer.
\epf
Au cours de la d\'emonstration, on a presque v\'erifi\'e
que la famille $\cal{C}$ a les propri\'et\'es suivantes :
\ble
\label{C6-4}
Les \'el\'ements de la famille alg\'ebrique $\cal{C}$ sont des
courbes contenues dans la surface $S$. Deux courbes g\'en\'eriques 
de la famille se coupent en $(n-1)$ points et, par toute famille 
g\'en\'erique de $n$ points de $S$, il passe une et
une seule courbe de la famille $\cal{C}$.
\ele
\bpf
Il reste \`a d\'emontrer la deuxi\`eme partie de l'\'enonc\'e.
Il r\'esulte du fait que la vari\'et\'e $\cal{C}$ est de dimension
$n$ que, par $n$ points g\'en\'eriques de $S$, il passe au moins une
courbe de $\cal{C}$. Comme deux courbes g\'en\'eriques 
se coupent en $(n-1)$ points, il en passe une seule.
\epf

\ssct{Conclusion}

Pour conclure, il suffit d'appliquer le r\'esultat suivant :
\bt[Enriques]
Soit $S$ une surface alg\'ebrique et $\cal{C}$ un syst\`eme alg\'ebrique 
de courbes sur $S$, de dimension $n$. Si la courbe g\'en\'erique de $\cal{C}$ est irr\'eductible 
et si, par une famille g\'en\'erique de $n$ points de $S$, il passe une et une seule 
courbe de $\cal{C}$, le syst\`eme $\cal{C}$ est un syst\`eme lin\'eaire.
\et
C'est en fait un \'enonc\'e de g\'eom\'etrie birationnelle
et l'\'enonc\'e, tel qu'il est formul\'e est incorrect. Le 
syst\`eme lin\'eaire final peut ne co\"{\i}ncider avec le syst\`eme 
initial que modulo un ensemble de dimension $<n$.

\sk
Quoi qu'il en soit, l'\'enonc\'e signifie qu'on peut choisir $\cal{C}=\P^n$
comme espace des param\`etres et que, pour tout point g\'en\'erique 
$s\in S$, l'ensemble des $\xi\in \cal{C}$
tels que $s$ appartient \`a la courbe $\xi$
est un hyperplan de l'espace projectif $\cal{C}$, identifi\'e \`a $\P^n$.

\sk
Comme on a dit dans l'introduction \`a cette section,
en associant \`a tout $x\in \C_0^n$ 
l'\'el\'ement de $\P^n$ qui correspond \`a la courbe $C(x)$
dans l'identification $\cal{C}=\P^n$, 
on obtient un diff\'eomorphisme local 
$\C_0^n \rightarrow \P^n$ qui lin\'earise le tissu.

%% file: bol.bbl
\begin{thebibliography}{aa}

\bibitem{B1} W. Blaschke, , 
    {\em Abh. Math. Semin. Hamb. Univ.} {\bf 9} (1933), 291--298.


\bibitem{B2} W. Blaschke, {\"U}ber die Tangenten einer ebenen Kurve f{\"u}nfter Klasse, 
    {\em Abh. Math. Semin. Hamb. Univ.} {\bf 9} (1933), 313--317.


\bibitem{BB} W. Blaschke, G. Bol, {\em Geometrie der Gewebe},
   Die Grundlehren der Mathematischen Wissenschaften, vol. 49, 
    {\em J. Springer, Berlin}, 1938 .

\bibitem{Bo}  G. Bol, Fl\"{a}chengewebe im dreidimensionalen Raum, 
    {\em Abh. Math. Semin. Hamb. Univ.} {\bf 10} (1934), 119--133.


\bibitem{E}  F. Enriques, Una questione sulla linearit\`a dei sistemi di 
curve appartenenti ad una superficie algebrica, 
    {\em Rend. Acc. Lincei} (1893), 3--8. 
      

\bibitem{CG1} S.S. Chern, P. Griffiths, Abel's theorem and webs,
  {\em Jahresber. Deutsch. Math.-Verein.} {\bf 80} (1978), {\em no.} 1-2, 13-110.


\bibitem{CG2} S.S. Chern, P. Griffiths, Corrections and addenda to our paper :
``Abel's theorem and webs'',  {\em Jahresber. Deutsch. Math.-Verein.} {\bf 83} (1981), {\em no.} 2, 78--83.

\bibitem{G1} P. Griffiths, Variations on a theorem of Abel, 
{\em Inventiones math.} {\bf 35} (1976), 321--390.


\bibitem{PT} L. Pirio, J.-M. Tr\'epreau, Tissus plans exceptionnels et fonctions th\^eta,
{\em Annales de l'Inst. Fourier}, \`a para\^{\i}tre.

\bibitem{P} H. Poincar\'e, Sur les surfaces de translation et les fonctions 
ab\'eliennes, {\em Bull. Soc. Math. France,} {\bf 29} (1901), 61--86.


\end{thebibliography}
